\documentclass[11pt]{article}
\usepackage[latin1]{inputenc}
\usepackage[english]{babel}
\usepackage{hyperref}
\usepackage{amsmath,amsfonts,amscd,amssymb}
\usepackage[all]{xy}
\usepackage{graphicx,epsfig}
\usepackage{makeidx}
\usepackage{color}
\usepackage[active]{srcltx} 

\newtheorem{thm}{Theorem}[section]
\newtheorem{rek}[thm]{Remark}

\newtheorem{lem}[thm]{Lemma}
\newtheorem{prop}[thm]{Proposition}
\newtheorem{cor}[thm]{Corollary}
\newtheorem{df}[thm]{Definition}

\newtheorem{probleme}[thm]{Problem}
\let\noi=\noindent

\def\notin{\mbox{$\in$ \hspace{-.8em}/}} 

\newcommand{\Cc}{\mathbb{C}}

\renewcommand {\epsilon}{\varepsilon}
\renewcommand {\le}{\leqslant}
\renewcommand {\ge}{\geqslant}
\renewcommand {\leq}{\leqslant}
\renewcommand {\geq}{\geqslant}

\renewcommand{\thefigure}{\arabic{figure}}

\makeatletter
\renewcommand{\fnum@figure}{{\figurename~\thefigure\ }}
\makeatother

\def\figurename{{Fig.}}%

\makeindex
\title{The Nash problem of arcs and the rational double point $\mathbf{E_6}$}
\author{\ by Camille Pl\'enat and Mark Spivakovsky}

\begin{document}

\maketitle

\begin{abstract}{This paper deals with the Nash problem, which consists in proving
that the number of families of arcs on a singular germ of a surface $S$
coincides with the number of irreducible components of the exceptional divisor in the minimal
resolution of this singularity. We propose a program for an affirmative solution of the
Nash problem in the case of normal 2-dimensional hypersurface singularities. We
illustrate this program by giving an affirmative solution of the Nash problem for
the rational double point $\mathbf{E_6}$. We also prove some results on the algebraic structure of the space of $k$-jets of an arbitrary hypersurface singularity and apply them to the specific case of $\mathbf{E_6}$.}
\end{abstract}


\thispagestyle{empty}

\thispagestyle{empty}

\pagenumbering{roman} \setcounter{page}{1}


\pagenumbering{arabic} \setcounter{page}{1}


\section{Introduction}
In this paper, $\Bbbk$ is an algebraically closed field of characteristic $0$.

Let $(S,0)$ be a normal surface singularity over $\Bbbk$ and $\pi:(X,E)\longrightarrow (S,0)$ the {\it minimal resolution} of
$(S,0)$, where $X$ is a smooth surface and $E=\pi^{-1}(0)$ is the exceptional set. Let $E=\bigcup\limits_{i\in \Delta}E_i$ be the decomposition of $E$ into its irreducible components, which we will
call {\it exceptional divisors}.

In order to study such a resolution, J. Nash (around 1968, published in 1995 \cite{N68})
introduced the space $H$ of {\it arcs} passing through the singular point $0$.
\begin{df} An \textbf{arc} is a $\Bbbk$-morphism from the local
ring $\mathcal{O}_{S,0}$ to the formal power series ring $\Bbbk[[t]]$.
\end{df}
Intuitively, an arc should be thought of as a parametrized formal curve, contained in $S$ and passing through the singular
point 0.

Nash had shown that $H$ has finitely many irreducible components, called {\it families of arcs}, and
that there exists a natural injective map, now called the Nash map, from the set of families of arcs to the set of
exceptional divisors of the minimal resolution. The celebrated Nash problem, posed in  \cite{N68}, is the question of
whether the Nash map is surjective.

Later on, M. Lejeune-Jalabert \cite{LEJ3} proposed the following decomposition of the
space $H$: let $N_i $ be the set of arcs whose strict transform in $X$ intersects
$E_i$ transversally but does not intersect any other exceptional
divisor $E_j$. M. Lejeune-Jalabert showed that $H=\bigcup\limits_{i\in \Delta}\overline{N_i}$
and the set $\overline{N_i}$ is an irreducible algebraic subvariety of the space of arcs;
therefore the families of arcs are among the $\overline{N_i}$'s. Moreover, notice that
there are as many $\overline{N_i}$ as divisors $E_i$. Then the Nash problem reduces
to showing that the $\overline{N_i}$ are precisely the irreducible components of $H$, that is, to proving
$card(\Delta)(card(\Delta)-1)$ non-inclusions:\\
\begin{probleme}\label{noninclusion}
Is it true that $\overline{N_i} \not \subset \overline{N_j}$ for all $i \not = j$?
\end{probleme}
This question has been answered affirmatively in the following special cases: for
$A_n$ singularities by Nash, for minimal surface singularities by A. Reguera
\cite{REG1} (with other proofs in J. Fernandez-Sanchez \cite{FER} and C.
Pl\'enat \cite{PLE1}), for sandwiched singularities by M. Lejeune-Jalabert and A.
 Reguera (cf. \cite{LEJ-REG} and \cite{REG2}), for toric vareties  by S. Ishii and J. Kollar
(\cite{ISH-KOL} using earlier work of C. Bouvier and G.
Gonzalez-Sprinberg \cite{BOU} and \cite{BOU-GON}), for rational double points
$\mathbf{D_n}$
by Pl\'enat \cite{PLE3}, for a family of non-rational surface singularities, as
well as for a family of singularities in dimension higher than $2$ by P.
Popescu-Pampu and C.Plénat  (\cite{P-PP1}, \cite{P-PP2}).\\In \cite{ISH-KOL}, S. Ishii and
J. Kollar gave a counter-example to the Nash problem in dimension greater than or
equal to $4$.

In this paper we prove the following theorem:

\begin{thm}
The Nash problem has an affirmative answer for rational double
points $\mathbf{E_6}$.
\end{thm}

But the principal aim of this paper is to present a general strategy
for attacking normal 2-dimensional hypersurface singularities which
has so far been successful in the case of $\mathbf{D_n}$ (\cite{PLE2})
and $\mathbf{E_6}$ (the present paper).

Once this theorem is proved, we have the following corollary
(cf. \cite{PLE1} for a proof):
\begin{cor}
Let $(S,0)$ be a normal surface singularity whose dual graph is
obtained from $\mathbf{E_6}$ by increasing the weights (that is,
allowing the exceptional curves to have self-intersection numbers of
the form $-n$ for $n\ge2$).
Then the problem also has an affirmative answer for $(S,0)$.
\end{cor}
Our program for solving the Nash problem for a normal 2-dimensional
hypersurface singularity with equation $F=\sum c_{\alpha \beta \gamma}
x^\alpha y^\beta z^\gamma =0$ is divided into two main steps.
For the first step we use the following valuative criterion:
\begin{prop}\label{critere}
Let $(S,0)$ be a normal surface singularity.\\
If there exists an element $f$ in $\mathcal{O}_{S,0}$ such that
$ord_{E_i}f<ord_{E_j}f$ then $\overline{N_i}\not \subset
\overline{N_j}$ .
\end{prop}
This result is stated and proved in (\cite{PLE3}, Proposition 1.1) for
arbitrary singularities in any dimension. It was first proved by
A. Reguera \cite{REG1} in a different, but equivalent formulation for
rational surface singularities.
\begin{rek}\label{Artintthm} Proposition \ref{critere} allows us to
  prove at least half of the non-inclusions appearing in Problem
  \ref{noninclusion} in the case of rational singularities. Indeed,
  let $(S,0)$ be a rational surface singularity and $E_i$, $E_j$ two
  distinct irreducible exceptional curves on the minimal resolution
  $X$ of $S$. Let $n=\#\Delta$. Since the intersection matrix
  $(E_q.E_s)$ is negative definite, there exists a cycle on $X$ of the
  form $C=\sum\limits_{q\in\Delta}m_qE_q$ such that
\begin{equation}
m_q>0,\ C.E_q\le0\quad\text{ for all }q\in \Delta\label{eq:fundamental}
\end{equation}
In fact, $n$-tuples $(m_1,\dots,m_n)$ of rational
numbers satisfying (\ref{eq:fundamental}) form an
$n$-dimensional cone in $\mathbb Q^n$, called the
Lipman cone. There exists a vector in the Lipman cone with integer
coefficients such that $m_i\ne m_j$, otherwise the Lipman cone would
be contained in the $(n-1)$-dimensional hyperplane $n_i=n_j$.  
Say, $m_i<m_j$. Since 
  $(S,0)$ is rational, Artin's theorem \cite{ART} tells us that there
  exists $f\in \mathcal O_{S,0}$ with $ord_{E_i}f=m_i$ and
  $ord_{E_j}f=m_j $, so the non-inclusion $\overline
  N_i\not\subset\overline N_j$ is given by the valuative
  criterion. This proves that for any pair $i,j\in \Delta$, $i\ne j$,
  at least one of the two non-inclusions $\overline
  N_i\not\subset\overline N_j$, $\overline N_j\not\subset\overline
  N_i$ is given by the valuative criterion.
\end{rek}
\bigskip

The second step consists in proving the remaining non-inclusions. For
this, we use the algebraic machinery developed in \S\ref{kjets} of
this paper. The idea is the following:\\
Let $E_i$ and $E_j$ be two exceptional divisors such that
\begin{equation}
ord_{E_i}f \leq ord_{E_j}f\text{  for\ all\ }f \in
\mathcal{O}_{S,0}.\label{eq:inequality}
\end{equation}
For rational surface singularities, the negative definiteness of the
intersection matrix $(E_i.E_j)$ implies that {\it strict} inequality
holds for at least one $f\in m_{S,0}$, so $\overline{N_i}\not \subset
\overline{N_j}$ by the valuative criterion (Proposition
\ref{critere}).\\
The opposite non-inclusion
\begin{equation}
\overline{N_j} \not\subset \overline{N_i}\label{eq:noninclusion}
\end{equation}
cannot be obtained from the valuative criterion and must be proved
separately.

Assume that $(S,0)$ is a normal hypersurface singularity, embedded in
the three-dimensional affine space $spec\ \Bbbk[x,y,z]$. An arc on
$(S,0)$ is described by three formal power series
\begin{equation}
\left\{\begin{array}{c}
\ x(t)=\sum\limits_{k=1}^\infty a_kt^k\\
\ y(t)=\sum\limits_{k=1}^\infty b_kt^k\\
\ z(t)=\sum\limits_{k=1}^\infty c_kt^k
\end{array}
\right.\label{eq:expansions}
\end{equation}
whose coefficients $a_k$, $b_k$, $c_k$ satisfy infinitely many
polynomial equations, obtained as follows. Substitute the series
(\ref{eq:expansions}) in $F$ and write
$F(x(t),y(t),z(t))=\sum\limits_{l=1}^\infty f_l(a,b,c)t^k$. Here
$a=(a_k)_{k\in\mathbb N}$, $b=(b_k)_{k\in\mathbb N}$,
$c=(c_k)_{k\in\mathbb N}$, and the $f_l$ are polynomials in $a$, $b$
and $c$. Let $\Bbbk^{\{a,b,c\}}$ denote the direct product of
infinitely many copies of $\Bbbk$, indexed by $a=(a_k)_{k\in\mathbb
  N}$, $b=(b_k)_{k\in\mathbb N}$ and $c=(c_k)_{k\in\mathbb N}$. We
think of $\Bbbk^{\{a,b,c\}}$ as an infinite-dimensional space over
$\Bbbk$ with coordinates $a,b,c$. Then $H$ is defined inside
$\Bbbk^{\{a,b,c\}}$ by the equations $f_l=0$, $l\in\mathbb N$.

To each arc as above we can associate in a natural way a closed point of the infinite-dimensional scheme $\mathcal
 H=Spec\frac{\Bbbk[a,b,c]}{(f)}$, where $(f)=(f_l)_{l\in\mathbb N}$. This scheme has the following description as a
 projective limit of schemes of finite type.\\
\begin{df}
An $i$-jet is a  $\Bbbk$-morphism $\mathcal{O}_{S,0}\rightarrow \frac{\Bbbk[[t]]}{(t^{i+1})}$.
\end{df}
Let us denote the set of all $i$-jets by $H(i)$.
 The set $H(i)$ can be naturally identified with the set of closed points of a scheme of finite type, denoted by $\mathcal
 H(i)$. With the natural maps $\rho_{ij}:\mathcal H(i)\rightarrow\mathcal H(j)$, $j<i$, called truncation maps, the
 $\mathcal H(i)$ form a projective system whose inverse limit is $\mathcal H$. The natural maps
 $\rho_i:\mathcal H\rightarrow\mathcal H(i)$ are also called \textbf{truncation maps}.

For a natural number $k$ and $i\in\Delta$, let $N_i(k)$ denote the image of $N_i$ in the algebraic variety $H(k)$ of $k$-jets of $S$.

We prove the non-inclusion (\ref{eq:noninclusion}) by contradiction:  suppose that
\begin{equation}
\overline{N_j} \subset \overline{N_i}.\label{inclusionienj}
\end{equation}
Clearly the inclusion (\ref{inclusionienj}) implies that $\overline{N_j(k)} \subset \overline{N_i(k)}$. Therefore we
 may work with $\mathcal H(k)$ for a sufficiently large $k$ instead of $\mathcal H$. The precise meaning of ``sufficiently
 large'' depends on the specific singularity in question, as well as on the particular non-inclusion (\ref{eq:noninclusion})
 we want to show; below we will specify $k$ precisely in each case.
Note that if $(S,0)$ is singular then $\rho_k$ need not, in general, be surjective onto $H(k)$.

Let $K(N_j(k))$ denote the field of rational functions of $N_j(k)$.

 By the Curve Selection Lemma (Lemma \ref{curveselection} below) there
exists a finite extension $L$ of $K(N_j(k))$ and an $L$-wedge
\begin{equation}
\phi_{ij}:Spec\frac{L[[t,s]]}{(t^{k+1})}\rightarrow S\label{eq:Lwedge}
\end{equation}
such that the image of the special arc $\{s=0\}$ is the generic point
of $N_j(k)$, while the image of general arc $\{s\ne0\}$ is an
$L$-point of $N_i(k)\setminus N_j(k)$. For each pair $i,j$ such that the non-inclusion (\ref{eq:noninclusion}) does not follow from the valuative criterion we study equations satisfied by an $L$-wedge (\ref{eq:Lwedge}) and prove that such an $L$-wedge does not exist.

The paper is organized as follows: in \S\ref{E6criterion} we first
recall the description of the singularity $\mathbf{E_6}$ we will use
and carry out the first step of the proof using the valuative
criterion. In \S\ref{kjets}, we partially describe the spaces of
$k$-jets $H(k)$ of a hypersurface singularity for a general $k$ and
apply this description to the specific case of the $\mathbf{E_6}$
singularity. We also describe the image of a family of arcs in the
truncated space $H(k)$. The last section is devoted to the second step
of the proof. Namely, we go one by one through the various
non-inclusions (\ref{eq:noninclusion}) which are not covered by the
valuative criterion and prove the non-existence of the $L$-wedge
(\ref{eq:Lwedge}) as above in each case. On four occasions, when the
resulting system of equations is too complicated to solve by hand, we
use MAPLE to check that it has no non-trivial solutions.

Note that by passing to the $k$-truncation we avoid using A. Reguera's
non-trivial theorem \cite{REG2}, which can be viewed as a version of
the Curve Selection Lemma for the pair of infinite dimensional schemes
$(N_i,N_j)$. In the present paper, the usual Curve Selection Lemma for
finite-dimensional algebraic varieties suffices for our purposes.

Recently several mathematicians announced positive solutions of the
Nash problem for more general types of singularities: Ana Reguera
for all rational surface singularities and Maria Pe Pereira (based on
the work \cite{BOB} of Javier Fernandez de Bobadilla) for quotients of
$\mathbb C^2$ by an action of finite group, though at the moment of
the writing of this manuscript their proofs have not yet been made
public. In any case, all the methods are completely different. We hope
that our method will one day be useful in a more general context, not
covered by the above results, such as normal hypersurface
singularities in $\mathbb C^3$.

\section{The singularity $\mathbf{E_6}$ and the valuative
  criterion}\label{E6criterion}

The singularity $\mathbf{E_6}$ is, by definition, the hypersurface
singularity defined in $\Bbbk^3$ by the equation  $F=z^2+y^3+x^4=0$.

The first graph in the following diagram is the dual graph of
$\mathbf{E_6}$; the remaining five graphs show the orders of vanishing
of the functions $x$, $y$, $z$, $z-ix^2$ and $z+ix^2$ on the
exceptional curves $E_1$, $E_2$, $E_3$, $E_4$, $E_5$ and $E_6$.

\begin{figure}[h,c]
$$
\setlength{\unitlength}{0.75cm}
\begin{picture}(14,5)
\includegraphics{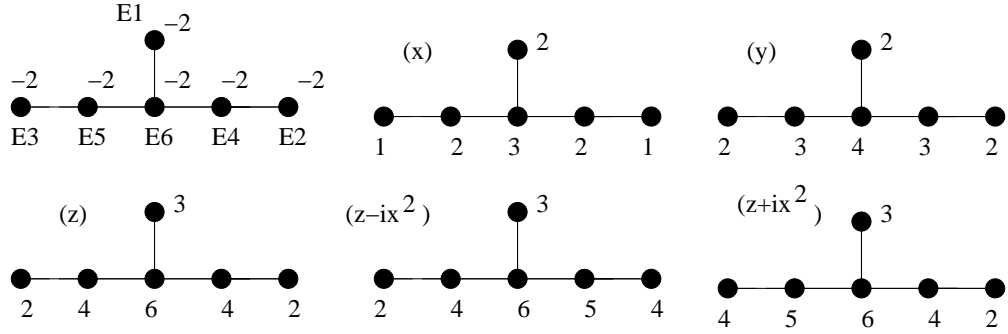}
\end{picture}
$$
\caption{ Dual graph of $\mathbf{E_6}$ and order of the functions $x$,
  $y$, $z$, $z-ix^2$ and $z+ix^2$}
\end{figure}
Consider the following partial ordering on the set
$\{E_1,E_2,E_3,E_4,E_5,E_6\}$. We say that
\begin{equation}
E_i<E_j\label{eq:Ei<Ej}
\end{equation}
if for all $f\in m_{S,0}$ the inequality (\ref{eq:inequality}) holds (as explained in Remark \ref{Artintthm}, together with the
rationality of $\mathbf{E_6}$ this implies that \textit{strict} inequality holds
in (\ref{eq:inequality}) for some $f\in m_{S,0}$).
 Using the functions $x$, $y$, $z$, $z-ix^2$ and $z+ix^2$, we see that our partial ordering contains at most the inequalities, shown in Figure 2.

\begin{figure}[ht]
$$
\setlength{\unitlength}{0.75cm}
\begin{picture}(6,4)
\includegraphics{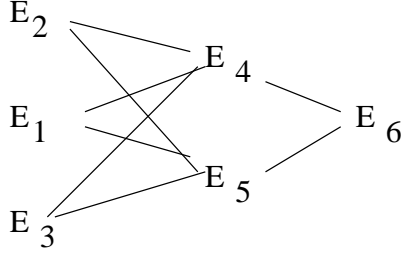}
\end{picture}
$$
\caption{The partial order for $\mathbf{E_6}$.}
\end{figure}

Here an inequality (\ref{eq:Ei<Ej}) is represented by placing $E_i$ to the left of $E_j$ (in fact, Figure 2 shows the entire partial ordering; this is all the information we can derive from comparing $ord_{E_i}f$ with $ord_{E_j}f$ for various $f \in
m_{S,0}$). Apply the valuative criterion with the functions $x$, $y$, $z$, $z-ix^2$ and $z+ix^2$. This proves all the non-inclusions (\ref{eq:noninclusion}) such that either
$E_j<E_i$ or $E_i$ and $E_j$ are not comparable in the partial ordering.
By symmetry, to complete the solution of the Nash problem for $E_6$, it is sufficient to show the following non-inclusions:
 \begin{eqnarray}
\overline{N_4}, \overline{N_6} &\not \subset& \overline{N_1}\label{eq:461}\\
\overline{N_4}, \overline{N_5},\overline{N_6} &\not \subset& \overline{N_2}\label{eq:4562}\\
\overline{N_6}&\not \subset& \overline{N_4}\label{eq:64}
\end{eqnarray}

For these non-inclusions we work in the space of $k$-jets  of the singularity $\mathbf{E_6}$ (with $k$ depending on the
non-inclusion). Let $\mathcal{P}_i$ and $\mathcal{P}_j$ be two prime ideals such that
\begin{eqnarray}
\overline{N_i(k))}&=&V(\mathcal{P}_i)\qquad\text{ and}\\
\overline{N_j(k))}&=&V(\mathcal{P}_j).
 \end{eqnarray}
 In order to prove that $N_i \not \subset N_j$, we show that $\mathcal{P}_j \not \subset \mathcal{P}_i$. To do this, we
 partially describe the ideals $\mathcal{P}_j$ and $\mathcal{P}_i$ and the space $H(k)$ of $k$-jets. This is the aim of
 \S\ref{kjets}.


 \section{The space of k-jets of a hypersurface singularity}\label{kjets}

 In this section we first recall some lemmas about hypersurface singularities, found in \cite{PLE3}. We then study the image
 of a family of arcs in the truncation spaces $H(k)$.\\
 In what follows we will look at a hypersuface singularity defined by
 $$
 f(x,y,z)=\sum{c_{\alpha \beta
\gamma}x^{\alpha}y^{\beta}z^{\gamma}}=0,
$$
embedded in $\Bbbk^3$ with a singularity at $0$.

\subsection{The $k$-jets scheme}

Fix an integer $k>0$.

   Any $k$-jet $\phi(t)$ passing through the singularity can be represented by
  three polynomials of degree $k$,
$\phi(t)=(x(t),y(t),z(t))=(a_1t+...+a_kt^k,b_1t+...+b_kt^k,c_1t+...+c_kt^k)$, with
  $$
  a_0=b_0 =c_0=0
  $$
  (because the singularity is at $0$), satisfying the following algebraic constraints.
Let $R_k $ be the polynomial ring $\Bbbk[a_1,...,a_k,b_1,...,b_k,c_1,...,c_k]$. Write
 $$
 (f\circ\phi)(t)=f(x(t),y(t),z(t))=0
 $$
  Let $f_i$ be the coefficient of $t^i$ in $f(x(t),y(t)z(t))$). Then $\{f_1=0,..,f_k=0\}$  are the
equations defining the k-jet scheme ${H}(k)$ in $\Bbbk^{3k}$.

 Let $l,m,n$ be integers such that there exists an exceptional divisor $E$ with
 \begin{eqnarray}
 ord_{E}(x)&=&l\\
 ord_E(y)&=&m\\
 ord_E(z)&=&n.
 \end{eqnarray}
  Let $K$ be the subset of the $k$-jet scheme defined in $H(k)$
  by the ideal $(a_1,...,a_{l-1},b_1,...,b_{m-1}$, $c_1$, ..., $c_{n-1})$.\\
 Let $r$ the smallest integer such that $$f_r \not \in
  (a_1,...,a_{l-1},b_1,...,b_{m-1},c_1,...c_{n-1}).$$\\
 The subspace $K$ of $H(k)$ ($k>r$) is defined by $$(a_1,...,a_{l-1},b_1,...,b_{m-1},c_1,...c_{n-1},f_r,...f_k)$$
in $\Bbbk^{3k}$.

 Then one can write, for $i \ge0$,
 \begin{equation*}
 \begin{split}
 f_{r+i}=\left(\frac{\partial f_r}{\partial a_l}\right)a_{l+i}+\left(\frac{\partial f_r}{\partial b_m}\right)&b_{m+i}+\left(\frac{\partial f_r}{\partial c_n}\right)c_{n+i}\\
 &+S_{r+i}(a_l,...,a_{l+i-1},b_m,...,b_{m+i-1},c_n,...,c_{n+i-1}),
\end{split}
 \end{equation*}
 where $S_{r+i}$ is a polynomial (for a proof see \cite{PLE3}, \S4.2).

Let us recall the main lemma of \cite{PLE3}, \S1.3, used for the description of the image of a family of arcs in the space of $k$-jets:

\begin{lem}\label{distinguished}
Consider the polynomial ring $R=\Bbbk[y_1,...y_n,x_{21},...,x_{2m},...,x_{k1},...,x_{km}]$, where $y_1,...y_n,x_{21},...,x_{2m},...,x_{k1},...,x_{km}$ are independent variables. Let $f_1,...,f_k$ be a sequence of elements of the following form :
\begin{eqnarray*}
\ f_1=f_1(y_1,...,y_n)=g_1...g_s\\
\ f_2=a_1x_{21}+...+a_mx_{2m}+h_2(y_1,...,y_n)\\
\ f_3=a_1x_{31}+...+a_mx_{3m}+h_3(y_1,...,y_n,x_{21},...,x_{2m})\\
\ \vdots \\
\ f_k=a_1x_{k1}+...+a_mx_{km}+h_k(y_1,...,y_n,x_{21},...,x_{(k-1)m})\\
\end{eqnarray*}
with $g_1,...,g_s$ distinct irreducible polynomials and $a_1,...,a_m \in \Bbbk[y_1,...,y_n]$.\\
For a fixed $j$, $1 \leq j \leq s$, let $S_j \subset \{ a_1,...,a_m\}$ be the set of
 $a_l$ such that $a_{l} \not \in (g_j)$.\\
Let us denote $I=(f_1,...,f_k)$ .\\
If $S_j \not = \emptyset$, there exists a unique minimal prime ideal $\mathcal{P}_j$ of $I$ such that $g_j \in \mathcal{P}_j$ and $a_l \not \in \mathcal{P}_j$ for all $a_l \in S_j$.\\
Assume $S_j \not = \emptyset$ for all $j\in\{1,\dots,s\}$. Let $\mathcal{Q} $ be a minimal prime ideal of $I$ different from
$ \mathcal{P}_1,...,\mathcal{P}_s$; then
 $(a_1,...,a_m) \subset \mathcal{Q}$.\\
Let $g_i$ and $g_j$ be two irreducible factors of $f_1$. Then
$\mathcal{P}_i \not = \mathcal{P}_j$. Finally, we have
$ht(\mathcal{P}_j)=k$.
\end{lem}
 \begin{df}
We call the prime ideal $\mathcal{P}_j$ of the lemma the distinguished ideal of $I$, associated to $g_j$.
\end{df}

Lemma \ref{distinguished} says that there are exactly $s$ distinguished
ideals of $I$, one associated to each irreducible factor $g_j$.

 \subsection{Image of a family of arcs in H(k)}\label{imageofarcs}
$  $\medskip

 Assume that $f$ is irreducible.
 Let $N_i$ be the set of arcs determined by the exceptional divisor $E_i$, as defined in the introduction.
 \medskip

 \noi{\bf Notation :}
 \begin{itemize}
  \item For an element $g\in\mathcal{O}_{S,0}$, let $\mu_i(g)$ be the order of vanishing of $g\circ\pi$ on $E_i$.
\item Let $R_k=\Bbbk[a_1,\dots,a_k,b_1,\dots,b_k,c_1,\dots,c_k]$.
 \item For $i\in\Delta$, let $o_i=\min\{ \alpha \mu_i(x)+\beta \mu_i(y)+\gamma \mu_i(z)\ |\ c_{\alpha \beta \gamma}\not = 0\}$.
 \item For $i\in\Delta$, $j\in\mathbb N$, let
\begin{eqnarray*}
 &o_{ij}=\min\{ [(\alpha-1) \mu_i(x)+\beta \mu_i(y)+\gamma \mu_i(z)]+j, [\alpha \mu_i(x)+(\beta-1) \mu_i(y)+\gamma \mu_i(z)]+j,\\
 &[\alpha \mu_i(x)+\beta \mu_i(y)+(\gamma-1) \mu_i(z)]+j\ |\  c_{\alpha \beta \gamma}\not = 0  \}.
\end{eqnarray*}
 \item Let $f_{j}$ be the coefficient of $t^j$ in $f(x(t),y(t),z(t))=0$
 \item Let $f_{ij}$ denote the unique element of $\Bbbk\left[a_{\mu_i(x)},\dots,a_j,b_{\mu_i(y)},\dots,b_j,
     c_{\mu_i(z)},\dots,c_j\right]$ such that $f_{ij}\equiv f_j$ modulo the ideal $(a_1,...,a_{\mu_i(x)-1},b_1,...,b_{\mu_i(y)-1},c_1,...,c_{\mu_i(z)-1})$ (here we adopt the obvious convention that the list $a_{\mu_i(x)},\dots,a_j$ is considered empty whenever $\mu_i(x)>j$, and similarly for the $b$ and $c$ coefficients).
 \end{itemize}

 \begin{prop}\label{tildeI}
  $  $
 Take an integer $k> o_i$.\\
 Let
 $$
 I_{ik}=\left(a_1,...,a_{\mu_i(x)-1},b_1,...,b_{\mu_i(y)-1},
c_1,...,c_{\mu_i(z)-1},f_{io_i},...,f_{io_{ik}}\right)R_k.
 $$
and
$$
\tilde I_{ik}=\left(I+ \left(a_1,...,a_{\mu_i(x)-1},b_1,...,b_{\mu_i(y)-1},
c_1,...,c_{\mu_i(z)-1}\right)\right)\cap R_k.
$$
 Then
 \begin{equation}
 I_{ik} \subset\tilde I_{ik}.\label{eq:Iik}
 \end{equation}
For $d\in\left\{\frac{\partial f_{io_i}}{\partial
a_{\mu_i(x)}},\frac{\partial f_{io_i}}{\partial b_{\mu_i(y)}},\frac{\partial
f_{io_i}}{\partial c_{\mu_i(z)}}\right\}$, we have
\begin{equation}
 I_{ik}(R_k)_d=\tilde I_{ik}(R_k)_d.\label{eq:Iik=}
 \end{equation}
\end{prop}
\noindent\textit{Proof.} The inclusion (\ref{eq:Iik}) is obvious. To prove
(\ref{eq:Iik=}), first note that the left hand side is contained in the right
hand side by (\ref{eq:Iik}). Conversely, let
$d=\frac{\partial f_{io_i}}{\partial a_{\mu_i(x)}}$; the proof for the other
two possible choices of $d$ is exactly the same.
Take an element $g\in\left(\tilde I_{ik}\right)_d$. By definition of $\tilde I_{ik}$, $g$
can be written in the form
 \begin{equation}
 g=\sum\limits_{j=o_i}^s h_jf_{ij}+\tilde g,\label{eq:s<}
 \end{equation}
where $h_j\in R_d$ and $\tilde
g\in\left(a_1,...,a_{\mu_i(x)-1},b_1,...,b_{\mu_i(y)-1},c_1,...,c_{\mu_i(z)-1}
\right)R_d$. Up to multiplication by a unit of $R_d$ (namely,
by $\frac1{\partial f_{io_i}/\partial a_{\mu_i(x)}}$), $f_{ij}$ has the form
$a_j+\lambda_{ij}$, where
$$
\lambda_{ij}\in\Bbbk[a_1,\dots,a_{j-1},b_1,\dots,b_j,c_1,\dots,c_j]_d
$$
Thus by adding a suitable multiple of $f_{ij}$ to each $h_{j'}$ with $j'<j$, we may
assume that $h_{j'}$ does not involve the variable $a_j$ whenever $j'<j$.
Also, we may assume that $\tilde g=0$ and that none of the $h_j$ involve the
variables $a_1,...,a_{\mu_i(x)-1},b_1,...,b_{\mu_i(y)-1},c_1,...,c_{\mu_i(z)-1}$.
We will now show that under these assumptions $s\le o_{ik}$ in (\ref{eq:s<}).
Indeed,  the right hand side of (\ref{eq:s<}) contains
exactly one term involving $a_s$. If we had $s>o_{ik}$ then, by
definition of $o_{ik}$, we have $g\notin R_k$, a contradiction.
This proves the equality (\ref{eq:Iik=}).\\
$\Box$\medskip

Let $\tau=\{\mu(x),\mu(y),\mu(z)\}$ be a triple such that there
exists $i\in \Delta$ with $$ \tau=\{\mu_i(x),\mu_i(y),\mu_i(z)\}. $$
Let $E(\tau)=\{E_l : \{\mu_l(x),\mu_l(y),\mu_l(z)\}=\tau\}$. For $E_i\in
E(\tau)$ and $j\in\mathbb N$, the numbers $o_i$, $o_{ij}$,
$\mu_i(x),\mu_i(y),\mu_i(z)$, the polynomials $f_{ij}$ and
the ideals $I_{ik}$, $\tilde I_{ik}$ depend only on $\tau$ and not on the particular
choice of $E_i\in E(\tau)$. We will therefore denote these objects by
$o_\tau$, $o_{\tau j}$, $\mu_\tau(x),\mu_\tau(y),\mu_\tau(z)$, $f_{\tau j}$, and $I_{\tau k}$, $\tilde I_{\tau k}$,
respectively.
\begin{prop}\label{image}{\bf (Image of a family)}
Assume that $f_{\tau o_\tau}$ is reduced but not necessarily irreducible and that it is not divisible by any of
$a_{\mu_\tau(x)},b_{\mu_\tau(y)},c_{\mu_\tau(z)}$; let $f_{\tau o_\tau}=g_1...g_s$ be its factorization into
irreducible factors.\\
Then:
\begin{itemize}
\item there exists a uniquely determined injective map
$$
\psi: \{1,...,s\}\longrightarrow E(\tau)
$$
such that for $j\in\{1,...,s\}$ and $E_i =\psi(j)$, the variety $\overline{N_i(k)}$ is defined by the distinguished prime
 ideal of $I_{\tau k}$  associated with $g_j$.
 \item  The non-inclusion (\ref{eq:noninclusion}) holds for all $E_i,E_j\in Im(\psi)$. In particular, if the map
 $\psi$ is surjective, (\ref{eq:noninclusion}) holds for all $E_i,E_j\in E(\tau)$.
\end{itemize}
\end{prop}
\begin{rek} If $s= card(E(\tau))$ then $\psi$ is necessarily bijective. This is the case for rational double points
$\mathbf{A_n}$, $\mathbf{D_n}$ (in both cases $s=card(E(\tau))=1$ for all values of $\tau$ \cite{PLE3}). Below, we will see
that for the singularity $\mathbf{E_6}$ we always have $s=card(E(\tau))\le2$, so, again, $\psi$ is bijective. Of course,
$\psi$ is bijective for any singularity for which the Nash problem has an affirmative answer. At this point we do not know
how to prove the bijectivity of $\psi$ for an arbitrary isolated 2-dimensional hypersurface singularity.
\end{rek}
\noi\textit{Proof of Proposition \ref{image}.} For the first assertion, note that the ideal $I_{\tau k}$ satisfies the
hypotheses of Lemma \ref{distinguished}, with the partial derivatives
$\frac{\partial f_{\tau o_\tau}}{\partial a_{\mu_\tau(x)}},\frac{\partial
f_{\tau o_\tau}}{\partial b_{\mu_\tau(y)}},\frac{\partial f_{\tau o_\tau}}{\partial
c_{\mu_\tau(z)}}$ playing the roles of $a_1,a_2,a_3$.

By definitions
\begin{equation}\label{eq:VtildeI}
V\left(\tilde I_{\tau k}\right)=\bigcup\limits_{\begin{array}{c}\mu_i(x)\ge\mu_\tau(x)\\ \mu_i(x)\ge\mu_\tau(x)\\
\mu_i(x)\ge\mu_\tau(x)\end{array}}\overline{N_i(k)}.
\end{equation}
Let $d$ be one of the partial derivatives of $f_{\tau o_\tau}$, which is not identically zero.
The fact that $f_{\tau o_\tau}$ is reduced implies that $I_{\tau k}R_d$ is not the unit
ideal. Now Proposition \ref{tildeI} (particularly, (\ref{eq:Iik=})) implies that the distinguished prime ideals
$\mathcal{P}_{jk}$, $j\in\{1,\dots,s\}$ of
$I_{\tau k}$ are also minimal primes of $\tilde I_{\tau k}$. Since the varieties $\overline{N_i(k)}$ are irreducible,
(\ref{eq:VtildeI}) shows that for each $j\in\{1,\dots,s\}$ there exists $i$ with
\begin{eqnarray}
\mu_i(x)&\ge&\mu_\tau(x),\\
\mu_i(x)&\ge&\mu_\tau(x),\\
\mu_i(x)&\ge&\mu_\tau(x),
\end{eqnarray}
such that $V(\mathcal{P}_{jk})=\overline{N_i(k)}$. Furthermore, since $g_j$ is
not divisible by $a_{\mu_\tau(x)}$, $b_{\mu_\tau(y)}$ or $c_{\mu_\tau(z)}$ and has no common factors with $d$ by assumption,
by Nullstellensatz there exist triples $(a',b',c')\in \Bbbk^3$ such that $g_j(a',b',c')=0$, $d(a',b',c') \not =0$ and $a',b',c'$ are different from
$0$. Then there exists an arc in $V(\mathcal{P}_{jk})$ of the form
 $\phi(t)=(a't^{\mu_\tau(x)}+\dots,b't^{\mu_\tau(y)}+\dots,c't^{\mu_\tau(z)}+\dots)$. Namely, we construct such an arc by
describing the values of $a_{\mu_\tau(x)+r}$, $b_{\mu_\tau(y)+r}$ and $c_{\mu_\tau(z)+r}$. We put
$(a_{\mu_\tau(x)},b_{\mu_\tau(y)},c_{\mu_\tau(z)})=(a',b',c')$. Then, for each positive integer $r$, we let
$b_{\mu_\tau(y)+r}$ and $c_{\mu_\tau(z)+r}$ be arbitrary elements of $\Bbbk$ and set
 $$
 a_{\mu_\tau(x)+r}=-\frac{f_{\tau,o_\tau+r}-a_{\mu_\tau(x)+r}d}d.
 $$
 This proves that $E_i\in E(\tau)$. We define $E_i=\psi(j)$.

The injectivity of $\psi$ is obvious from the definition. Also by definition, the non-inclusion (\ref{eq:noninclusion}) is
satisfied for all $E_i,E_j\in Im(\psi)$. Thus, if $\psi$ is surjective, (\ref{eq:noninclusion}) holds for all $E_i,E_j\in
E(\tau)$, as desired. This completes the proof.
$\Box$\medskip

 \noindent\textbf{Example.} Let us apply the above ideas to the special case of the $\mathbf{E_6}$ singularity. According
 to Figure 1, there are four possible values of $\tau$: (2,2,3), (1,2,2), (2,3,4) and (3,4,6). We have $E(2,2,3)=\{E_1\}$,
$E(1,2,2)=\{E_2,E_3\}$, $E(2,3,4)=\{E_4,E_5\}$, and $E(3,4,6)=\{E_6\}$. Thus, for $\tau=(2,2,3)$ or $\tau=(3,4,6)$ the
bijectivity of the map $\psi$ is immediate.

Next, let $\tau=(1,2,2)$. We have $o_\tau=4$ and $f_{\tau o_\tau}=c_2^2+a_1^4=(c_2+ia_1^2)(c_2+ia_1^2)$, so
$f_{\tau o_\tau}$ is a product of two distinct irreducible factors.

Similarly, if $\tau=(2,3,4)$, we have $o_\tau=8$ and $f_{\tau o_\tau}=c_4^2+a_2^4=(c_4+ia_2^2)(c_4+ia_2^2)$, so, again
$f_{\tau o_\tau}$ is a product of two distinct irreducible factors.

Since in the last two cases $f_{\tau,o_\tau}$ has two irreducible factors and $\#E(\tau)=2$, the map $\psi$ is bijective
also in these two cases. It follows from Proposition \ref{image} that for a sufficiently large $k$ each $\overline{N_i(k)}$ is of
the form $V(\mathcal P_{ik})$, where $\mathcal P_{ik}$ is a distinguished prime ideal, associated to $I_{ik}$.

We recall that the goal is to prove that
\begin{equation}
\mathcal{P}_{ik} \not \subset \mathcal{P}_{jk}\label{eq:Pnoninclusion}
\end{equation}
whenever
\begin{equation}
E_i<E_j.\label{eq:EiEj}
\end{equation}

\subsection{The strategy for proving the non-inclusion (\ref{eq:Pnoninclusion})}\label{strategy}

$ $
\bigskip

By the valuative criterion we already have the opposite non-inclusion in
(\ref{eq:Pnoninclusion}). Inequality (\ref{eq:EiEj}) means that
$ord_{E_i} g \leq ord_{E_j}g$ for all $g \in \mathcal{O}_{S,0}$. We thus have the
following inclusions:
$$
\left\{ \begin{array}{c}
\ I_{ik} \subset \mathcal{P}_{ik}\\
\ I_{ik}\subset I_{jk} \subset \mathcal{P}_{jk}
\end{array} \right.
$$
\bigskip

Assume ${N_j}(k)  \subset  \overline{N_i(k)}$ for a certain order $k$.
\bigskip

We will need the Curve Selection lemma (for usual finite-dimensional
algebraic varieties):

\begin{prop}\label{curveselection}\textbf{(Curve Selection Lemma).}
Let $V$ be a reduced algebraic variety over an algebraically closed field $\Bbbk$ and $W$ a proper reduced irreducible
subvariety of $V$. Let $K(W)$ denote the field of rational functions
of $W$.\\
There exists a finite field extension $L$ of $K(W)$ and an arc $\phi:Spec\ L[[s]]\rightarrow V$ whose generic point
maps to $V\setminus W$, and the special point to the generic point of $W$.
\end{prop}
\noi\textit{Proof:} Replacing $V$ by a suitable affine open subset of it, we may assume, without loss of generality, that $V$ is an affine variety. Let $A$ denote the coordinate ring of $V$ and write $W=V(P)$ where $P$ is a prime ideal of $A$. Let $Q$ denote a prime ideal of $A$, contained in $P$, such that $ht\ Q=ht\ P-1$. Let $B$ denote the normalization of the ring $\frac{A_P}{QA_P}$, $\hat B$ the completion of $B$ at some fixed maximal ideal and $L$ the residue field of $\hat B$. The field $L$ is a finite extension of $K(W)$. Then $\hat B$ is a complete regular 1-dimensional local ring; let $s$ be a regular parameter of $\hat B$. We have $\hat B\cong L[[s]]$; the composition of the natural maps $A\rightarrow A_P\rightarrow\frac{A_P}{QA_P}\rightarrow B\rightarrow\hat B$ induces the morphism $\phi$ required in the Proposition.
$\Box$\medskip

Let $W=N_j(k)$. In our context, the curve is an arc of the form
$\phi_{ij}:Spec\ L[[s]]\rightarrow N_i(k)$, which corresponds to a
``truncated'' $L$-wedge
\begin{equation}
\phi_{ij}:Spec\frac{L[[t,s]]}{(t^{k+1})} \rightarrow (S,0)\label{eq:coin}
\end{equation}
with special arc ($s=0$) maps to the generic arc of $N_j(k)$ and whose general
arc maps to an $L$-point of $N_i(k)\setminus N_j(k)$. A wedge as in (\ref{eq:coin}) is given by three
polynomials of the form
\begin{eqnarray}
x(t,s)&=&\sum\limits_{n=0}^k\mathbf a_n(s)t^n\\
y(t,s)&=&\sum\limits_{n=0}^k\mathbf b_n(s)t^n\\
z(t,s)&=&\sum\limits_{n=0}^k\mathbf c_n(s)t^n
\end{eqnarray}
Write the coefficients $\mathbf a_n(s),\mathbf b_n(s),\mathbf c_n(s)$ of the wedge in the form
$$
\left\{\begin{array}{c}
\ \mathbf a_n(s)=\sum\limits_{p=0}^\infty a_{np}s^p\\
\ \mathbf b_m(s)=\sum\limits_{p=0}^\infty b_{mp}s^p\\
\ \mathbf c_l(s)=\sum\limits_{p=0}^\infty c_{lp}s^p,\\
\end{array}
\right.
$$
with $a_{np},b_{mp},c_{lp}\in L$, where $a_{n0}, b_{m0},c_{l0}$
satisfy the equations of $N_j(k)$. In particular, $a_{n0}=0$ when
$n<ord_j(x)$, $b_{m0}=0$ when $m<ord_j y$ and $c_{l0}=0$ when $l<ord_j
z$. Let us denote by $\alpha_n$ (resp. $\beta_m$ and $\gamma_l$) the
smallest order $q$ for which  $a_{nq}$ (resp. $b_{mq}$ and $c_{lq}$)
is not $0$. We need to compute these exponents in order to construct
the wedge $\phi_{ij}$. Note that $a_{n0}\ne0$ if and only if
$\alpha_n=0$, and similarly for the $b$ and $c$ coefficients; we
always have $a_{n0} \ne 0$ if $n=ord_j x$.

The morphism (\ref{eq:coin}) is given by a ring homomorphism
\begin{equation}
\mathcal O_{S,0}\rightarrow \frac{L[[t,s]]}{(t^{k+1})}.\label{eq:dualcoin}
\end{equation}
Localizing $\frac{L[[t,s]]}{(t^{k+1})}$ by the element $s$, we obtain
an $L((s))$-point of $N_i(k)$ (informally, an $L((s))$-arc lying in
$N_i(k)$). Thus the coefficients $\mathbf a_n(s),\mathbf b_m(s),\mathbf c_l(s)$ satisfy the
equations $f_{iu}$ of $N_i(k)$ and their constant terms $a_{n0},b_{m0},c_{l0}$ satisfy the
equations $f_{ju}$ of $N_j(k)$ (here $f_{iu}$ is the coefficient of $t^u$ in $Fo\phi_{ij}$ and similarly for $f_{ju}$; see the
beginning of \S\ref{imageofarcs} where this notation was
introduced).

Let $A_{np}$, $B_{mp}$, $C_{lp}$, $p\ge0$, be independent variables and write
$$
\left\{\begin{array}{c}
\ A_n(s)=\sum\limits_{p=0}^\infty A_{np}s^p\\
\ B_m(s)=\sum\limits_{p=0}^\infty B_{mp}s^p\\
\ C_l(s)=\sum\limits_{p=0}^\infty C_{lp}s^p.\\
\end{array}
\right.
$$
We have finitely many equalities of the form
 \begin{equation}
0=f_{iu}(A(s),B(s),C(s))=\sum_{v=0}^\infty f'_{iv u}s^v,\quad u\le o_{ij},
\end{equation}
where $A(s)$ stands for $\{A_n(s)\}_{n\in\mathbb N}$, and similarly for $B$ and $C$.
Here the coefficients $f'_{iv u}$ are polynomials in $A_{np}$, $B_{mp}$, $C_{lp}$ which vanish after substituting $A_{np}=a_{np}$, $B_{mp}=b_{mp}$, $C_{lp}=c_{lp}$.
\bigskip

Let $J$ denote the ideal of $L[A,B,C]$ generated by all the elements
of the form $A_{np}$ with $p<\alpha_n$, $B_{mp}$ with $p<\beta_m$ and
$C_{lp}$ with $p<\gamma_l$, where $A$ stands for
$\{A_{np}\}_{p\in\mathbb N}$, and similarly for $B$ and $C$. Let
$\theta_u=\min\{v\ \vert\ f'_{iv u}(A,B,C)\notin J\}$. Write
$g_{\theta_u}=f'_{i\theta_u u}$. In other words,  $g_{\theta_u}$ is
the first non-zero coefficient of $f_{i,u}(A(s),B(s),C(s))$, viewed as
a series in $s$, not belonging to the ideal $J$.
\medskip

\noi\textbf{Notation.} For the rest of this paper, we will write $a_n$
for $a_{n\alpha_n}$, $b_m$ for $b_{m\alpha_m}$ and $c_l$ for
$c_{l\gamma_l}$.
\medskip

\begin{rek}
\begin{itemize}
\item The coefficient $g_{\theta_u}$  depends only on $A_{n\alpha_n}$,
  $B_{m \beta_m}$ and $C_{l \gamma_l}$. Since $a_n\ne0$, $b_m\ne0$,
  $c_l\ne0$ and
  $$
  g_{\theta_u}(a_n, b_m,c_l)=0,
  $$
  the coefficient $g_{\theta_u}$ cannot be a monomial in
  $a_n$, $b_m$, $c_l$. In general, $g_{\theta_u}$ is a
  quasi-homogeneous polynomial in which $A_{n\alpha_n}$ has weight
  $\alpha_n$, $B_{m\beta_m}$ weight $\beta_m$ and $C_{l\gamma_l}$
  weight $\gamma_l$. Equality of weights of different monomials
  appearing in $g_{\theta_u}$ will give us a system of conditions on
  the exponents $\alpha_n$,  $\beta_m$ and $\gamma_l$. More precisely,
  we are not interested in the values of $\alpha_n$,  $\beta_m$ and
  $\gamma_l$ \textit{per se} but rather in the ratios of the form
  $\frac{\alpha_n}{\delta}$, where $\delta$ is some fixed element of
  the set
  $\{\alpha_{\mu_i(x)},\beta_{\mu_i(y)},\gamma_{\mu_i(z)}\}$. In other
  words, we are interested in the ``normalized" weights $\alpha_n$,
  $\beta_m$ and $\gamma_l$, where we set, for example, the first
  non-trivial weight $\alpha_{\mu_i(x)}$ equal to 1.
\item The hardest part of the proof is to recover the coefficients
  $g_{\theta_u}$. In order to do this, we will use the fact that
  $g_{\theta_u}$ are not monomials to give lower bounds on $\alpha_n$,
  $\beta_m$ and $\gamma_l$.
\end{itemize}
\end{rek}

 The equation $g_{\theta_u}=0$ plus the equations $f_{jk}(a_{n0},
 b_{m0}, c_{l0})=0$ form a system satisfied by the coefficients of the
 wedge. If this system has no solutions then the wedge does not
 exist. In one exceptional case, that of the non-inclusion
 $\overline{N_4} \not \subset  \overline{N_2}$, we will need to use
 $f'_{i,\theta_u+1,u}$, the next coefficient of
 $f_{iu}(A_n(s),B_m(s),C_l(s))$ after $g_{\theta_u}$, to arrive at a
 contradiction.

\bigskip

In the next section we compute the weights $\alpha_n$ , $\beta_m$ and
$\gamma_l$ for the singularity $E_6$ and show that the system
$$
\left\{\begin{array}{c}
\ g_{\theta_u}=0\\
\ f_{ju}(a_{n0}, b_{mu}, c_{l0})=0\\
\end{array}\right.
$$
for the remaining non-inclusions other than $\overline{N_4} \not \subset  \overline{N_2}$, as well as the augmented system
$$
\left\{\begin{array}{c}
\ g_{\theta_u}=0\\
\ f'_{i,\theta_u+1,u}=0\\
\ f_{ju}(a_{n0}, b_{m0}, c_{l0})=0\\
\end{array}\right.
$$
in the case of the non-inclusion $\overline{N_4} \not \subset  \overline{N_2}$, have no solutions.

\bigskip

\section{ Computations and proof for the $\mathbf{E_6}$ singularity}

Let us consider the $\mathbf{E_6}$ singularity and study the different non-inclusions.
For each non-inclusion $\overline{N_j}\not \subset \overline{N_i}$ appearing in
(\ref{eq:461})--(\ref{eq:64}), we will denote $R(k)=\frac{R_k}{\mathcal{P}_{ik}}$.

\noi{\bf Notation}: When talking about the non-inclusion $\overline{N_j}\not \subset \overline{N_i}$, the notation $a\ |\ b$ will mean ``a divides b in $\overline{R(k)}$", unless otherwise specified (here $\overline{R(k)}$ stands for the integral closure of $R(k)$ in its field of fractions). For some non-inclusions, we will study
divisibility in a suitable localization of $\overline{R(k)}$, which will be specified explicitly in each case.

For each of the six non-inclusions involved, it is sufficient to prove that
\begin{equation}
 \mathcal{P}_{ik} \not\subset \mathcal{P}_{jk}\label{eq:PiknotinPij}
\end{equation}
for  some $k$, in particular for $k=o(j)$. Take $k=o(j)$.

We prove the non-inclusion (\ref{eq:PiknotinPij}) by contradiction. Assume that
 $\mathcal{P}_{ik} \subset \mathcal{P}_{jk}$. By the Curve Selection lemma there exists an $L$-wedge whose special arc is the generic point of $N_j(k)$ and whose generic arc is in $N_i(k)$. The first coefficient $g_{\theta_u}$ of $f_{iu}$ cannot be a monomial as generically on $N_i(k)$ each monomial in $a_n$, $b_m$, $c_l$ is not zero.
\bigskip

As explained above, we are interested in computing ratios of the form
$\frac{\alpha_n}{\delta}$, where $\delta$ is some fixed element of the
set $\{\alpha_{\mu_i(x)},\beta_{\mu_i(y)},\gamma_{\mu_i(z)}\}$, and
$\mu_i(x)\le n<\mu_j(x)$, and similarly for $\frac{\beta_m}{\delta}$,
$\mu_i(y)\le m<\mu_j(y)$, and $\frac{\gamma_l}{\delta}$, $\mu_i(z)\le
l<\mu_j(z)$ (we will pick and fix a specific $\delta$ in the proof of
each non-inclusion, but the choice of $\delta$ will depend on the
non-inclusion we want to prove). For example, suppose
$\delta=\alpha_{\mu_i(x)}$. Then our problem is closely related to
studying, for each $n$, the totality of pairs
$(\alpha,\delta')\in\mathbb N^2$ such that
\begin{equation}
\mathbf a_n(s)^\alpha\ \left|\ \mathbf
  a_{\mu_i(x)}(s)^{\delta'}\right.,\label{eq:andivides}
\end{equation}
and similarly for $\mathbf b_m(s)^\beta\ \left|\ \mathbf
  a_{\mu_i(x)}(s)^{\delta'}\right.$ and $\mathbf c_l(s)^\alpha\
\left|\ \mathbf a_{\mu_i(x)}(s)^{\delta'}\right.$. Precisely, we have
$$
\frac{\alpha_n}{\alpha_{\mu_i(x)}}=
\inf\left\{\frac\alpha{\delta'}\right\},
$$
where $(\alpha,\delta')$ runs over all the pairs satisfying
(\ref{eq:andivides}).
\bigskip

\begin{rek} In \cite{PLE2} and \cite{PLE3} a different method is used
  to prove the non-inclusions not covered by the valuative
  criterion. Namely, we use the fact that the ideal $\mathcal P_{ik}$
  can be expressed as the saturation $(\mathcal P_{ik}R(k):d^\infty)$,
  where $d\in\{a_{\mu_i(x)},b_{\mu_i(y)},c_{\mu_i(z)}\}$. For most
  non-inclusions, we explicitly construct elements of $(\mathcal
  P_{ik}R(k):d^\infty)$, not belonging to $\mathcal P_{jk}$, which
  settles the problem. In both the saturation and the wedge methods,
  the key point is to compute the weight ratios of the form
  $\frac{\alpha_n}{\delta}$, $\frac{\beta_m}{\delta}$ and
  $\frac{\gamma_l}{\delta}$ as above. One advantage of the wedge
  method is that it gives a more geometric vision of the proof.
\end{rek}
\bigskip

In what follows we truncate at the order $o_j$.
\begin{enumerate}
\item $\bullet$ {\bf $\overline{N_4} \not \subset \overline{N_1}$}.
In this case we truncate at the order $o_4=8$. We have $o_1=6$.

Assume that $\overline{N_4(8)} \subset \overline{N_1(8)}$, aiming for
contradiction.
Let $\phi_{42}$ be a wedge with generic arc living in $N_1(8)$ and
special arc mapping to $N_4(8)$. Then the wedge is of the form:

\begin{equation}\label{eq:wedge}
\left\{\begin{array}{c}
\mathbf b_2(s)={b_2}s^{\beta_2}+\sum\limits_{q=\beta_2+1}^\infty b_{2q}s^q\\
\mathbf c_3(s)={c_3}s^{\gamma_3}+\sum\limits_{q=\gamma_3+1}^\infty c_{3q}s^q\\
\mathbf a_n(s)={a_n}+\sum\limits_{q=1}^\infty a_{nq}s^q,\quad n\ge2\\
\mathbf b_m(s)={b_m}+\sum\limits_{q=1}^\infty b_{mq}s^q,\quad m\ge3\\
\mathbf c_l(s)={c_l}+\sum\limits_{q=1}^\infty c_{lq}s^q,\quad l\ge4\\
\end{array}\right.
\end{equation}
where ${a_n}, {b_m}, {c_{l}}$ satisfy the equations of $N_4(6)$, and are non-zero elements of $L$.

The following equations hold on
$N_1(8)$:
\begin{eqnarray*}
\ \mathbf a_1=\mathbf b_1=\mathbf c_1=\mathbf c_2=0\\
\ f_{1,6}=\mathbf c_3^2+\mathbf b_2^3=0\\
\ f_{1,7}= 2\mathbf c_3\mathbf c_4+3\mathbf b_2^2\mathbf b_3=0.
\end{eqnarray*}

The following equations hold on $N_4(8)$:

\begin{eqnarray*}
\ {a_1}={b_1}={c_1}={c_2}={c_3}={b_{2}}=0\\
\ f_{4,8}={c_{4}}^2+{a_{2}}^4=0.\\
\end{eqnarray*}

The generic arc lives in $N_1(8)$, and thus satisfies the equations of $N_1(8)$. This leads to finitely many equations (as we are in $R(8)$):
 \begin{eqnarray*}
\ 0=f_{1,6}(\mathbf a(s),\mathbf b(s),\mathbf c(s))={c_3}^2s^{2\gamma_3}+{b_2}^3s^{3\beta_2}+\hdots\\
\ 0=f_{1,7}(\mathbf a(s),\mathbf b(s),\mathbf c(s))=2{c_3}{c_{4}}s^{\gamma_3}+3{b_2}^2{b_{3}}s^{2\beta_2}+\hdots\\
\vdots
\end{eqnarray*}
As $c_3 \not = 0$ and $b_2 \not = 0$, we obtain a relation between $\gamma_3$ and $\beta_2$ :
$$
2\gamma_3=3\beta_2
$$
which implies that
$$
\beta_2\le \gamma_3 <2\beta_2
$$

and that
\begin{eqnarray*}
\ {c_3}^2+{b_2}^3=0\\
\ 2c_3{c_{4}}=0.
\end{eqnarray*}
Thus for the equation $f_{4,8}$ we have $\theta_8=\gamma_3$ and $g_{\theta_8}=2c_3c_4$, which is impossible.


\bigskip

\item $\bullet$ {\bf $\overline{N_5} \not \subset \overline{N_2}$}.
In this case we truncate at the order $o_5=8$. We have $o_2=6$. Assume that $\overline{N_5}  \subset \overline{N_2}$, aiming for contradiction.
Let $\phi_{52}$ be a wedge  with generic arc living in $N_2(8)$ and special arc mapping to the generic arc in $N_5(8)$.

The following equations hold on
$N_2(8)$:
\begin{eqnarray*}
\ \mathbf b_1=\mathbf c_1=0\\
\ f_{2,4}=\mathbf c_2^2+\mathbf a_1^4=0 \\
\  f_{2,5}=2\mathbf c_2\mathbf c_3+4\mathbf a_1^3\mathbf a_2=0\\
\ f_{2,6}=\mathbf c_3^2+2\mathbf c_2\mathbf c_4+\mathbf b_2^3+4\mathbf a_1^3\mathbf a_3+6\mathbf a_1^2\mathbf a_2^2=0\\
\ f_{2,7}= 2\mathbf c_3\mathbf c_4+2\mathbf c_2\mathbf c_5+3\mathbf b_2^2\mathbf b_3+4\mathbf a_1^3\mathbf a_4+12 \mathbf a_1^2\mathbf a_2\mathbf a_3+4\mathbf a_2^3\mathbf a_1=0
\end{eqnarray*}
We have $f_{2,4}=(\mathbf c_2+i\mathbf a_1^2)(\mathbf c_2-i\mathbf a_1^2)$. As can be seen from Figure 1, $\mathcal P_{2,8}$ is the distinguished ideal
corresponding to the irreducible factor $\mathbf c_2-i\mathbf a_1^2$. Let us use the notation
$$
g_{2,2}:= \mathbf c_2-i\mathbf a_1^2.
$$
Combining $f_{2,5}$ and $g_{2,2}$ we see that $2i\mathbf c_3\mathbf a_1^2+4\mathbf a_1^2\mathbf a_2=0$ on $N_2(8)$. Since $\mathbf a_1$ does not vanish identically on $N_2(8)$, we have
$$
\bar f_{2,3}:=\mathbf c_3-2i\mathbf a_1\mathbf a_2=0
$$
on $N_2(8)$.


We claim that
\begin{eqnarray}
\gamma_2&=&2\alpha_1,\label{eq:gamma2}\\
\gamma_3&=&\alpha_1\label{eq:gamma3}\\
\beta_2&\ge&\frac23\alpha_1\label{eq:beta2}
\end{eqnarray}
Now,
\begin{itemize}
\item (\ref{eq:gamma2}) holds thanks to the equation $g_{2,2}=0$.
\item (\ref{eq:gamma3}) holds by the equation $\bar f_{2,3}=0$ and the fact that $\alpha_2=0$.
\item (\ref{eq:beta2}) holds by the equation $\bar f_{2,6}=0$, (\ref{eq:gamma2}) and (\ref{eq:gamma3}).
\end{itemize}
After a suitable automorphism of $L[[s]]$, we may assume that $a_1=1$.

The vanishing of the first non-trivial coefficients of the power series $\bar f_{2,3}(\mathbf a_1(s),\mathbf a_2(s),\mathbf c_(s))$ and $f_{2,7}$ gives the equations
\begin{eqnarray}
\  c_3-2ia_2=0\label{eq:c3intermsofa2}\\
\ 2c_3{c_4}+4{a_2}^3=0\label{eq:c4intermsofa2}
\end{eqnarray}
and we have (first equation of $N_5$):
\begin{equation}
{c_4}+i{a_2}^2=0.\label{eq:c4intermsofa2bis}
\end{equation}
Substituting (\ref{eq:c3intermsofa2}) into (\ref{eq:c4intermsofa2}) and dividing through by $4ia_2^2$, we obtain the equation
$$
{c_4}-i{a_2}^2=0,
$$
which contradicts (\ref{eq:c4intermsofa2bis}) and the fact that $c_4$ and $a_2$ are non-zero elements of $L$.

\item $\bullet$ {\bf $\overline{N_4} \not \subset  \overline{N_2}$. }
In this case we truncate at the order $o_4=8$.

Assume that ${N_4}(8)  \subset  \overline{N_2(8)}$, aiming for contradiction. We can construct  an  $L$-wedge $Spec\ L[[t,s]]\rightarrow E_6$, with the special arc mapping to the generic arc of $N_4$ and with the general arc lifting to $E_2$.

The following equations hold on $N_2(8)$:
\begin{eqnarray*}
\mathbf b_1=\mathbf c_1=0\\
g_{2,2}= \mathbf c_2-i\mathbf a_1^2=0 \\
\bar f_{2,3}=\mathbf c_3-2i\mathbf a_1\mathbf a_2=0\\
f_{2,6}=\mathbf c_3^2+2\mathbf c_2\mathbf c_4+\mathbf b_2^3+4\mathbf a_1^3\mathbf a_3+6\mathbf a_1^2\mathbf a_2^2=0\\
f_{2,7}= 2\mathbf c_3\mathbf c_4+2\mathbf c_2\mathbf c_5+3\mathbf b_2^2\mathbf b_3+4\mathbf a_1^3\mathbf a_4+12 \mathbf a_1^2\mathbf a_2\mathbf a_3+4\mathbf a_2^3\mathbf a_1=0\\
f_{2,8}=\mathbf c_4^2+2\mathbf c_3\mathbf c_5+2\mathbf c_2\mathbf c_6+3\mathbf b_2^2\mathbf b_4+\mathbf a_2^4+4\mathbf a_1^3\mathbf a_5+12\mathbf a_1^2\mathbf a_2\mathbf a_4+12\mathbf a_1\mathbf a_2^2\mathbf a_3+6\mathbf a_1^2\mathbf a_3^2=0.
\end{eqnarray*}
Modifying $f_{2,6}$ and $f_{2,7}$ by suitable multiples of $g_{2,2}$ and $\bar f_{2,3}$, we may replace them by
\begin{eqnarray*}
\bar f_{2,6}:=2i\mathbf a_1^2(\mathbf c_4-i\mathbf a_2^2)+\mathbf b_2^3+4\mathbf a_1^3\mathbf a_3=0\\
\bar f_{2,7}:=4i\mathbf a_1\mathbf a_2(\mathbf c_4-i\mathbf a_2^2)+2\mathbf c_2\mathbf c_5+3\mathbf b_2^2\mathbf b_3+4\mathbf a_1^3\mathbf a_4+12 \mathbf a_1^2\mathbf a_2\mathbf a_3=0\\
\bar f_{2,8}=\mathbf c_4^2+4i\mathbf a_1\mathbf a_2\mathbf c_5+2i\mathbf a_2\mathbf c_6+3\mathbf b_2^2\mathbf b_4+\mathbf a_2^4+4\mathbf a_1^3\mathbf a_5+12\mathbf a_1^2\mathbf a_2\mathbf a_4+12\mathbf a_1\mathbf a_2^2\mathbf a_3+6\mathbf a_1^2\mathbf a_3^2=0
\end{eqnarray*}
Note that the equation $f_{4,8}=c_4-ia_2^2=0$ vanishes on $\overline{N_4(8)}$.

Let $\mu$ denote the $s$-adic valuation of $L[[s]]$. We define $\alpha:= \mu(\mathbf c_4(s)-i\mathbf a_2(s)^2)$. We claim that
\begin{eqnarray}
\gamma_2&=&2\alpha_1\label{eq:gamma21}\\
\gamma_3&=&\alpha_1\label{eq:gamma31}\\
\beta_2&\ge&\alpha_1\label{eq:beta21}\\
\alpha&\ge&\alpha_1.\label{eq:alpha}
\end{eqnarray}
Indeed,
\begin{itemize}
\item (\ref{eq:gamma21}) holds thanks to the equation $g_{22}=0$.
\item (\ref{eq:gamma31}) is given by $\bar f_{2,3}$ and the fact that $a_{20}=a_2\ne0$, and hence
\begin{equation}
\alpha_2=0.\label{eq:a2invertible}
\end{equation}
\end{itemize}
We have $\alpha_1>0$. Using (\ref{eq:a2invertible}) once again, we obtain from the equations $\bar f_{2,6}=0$ and $\bar f_{2,7}=0$ that
\begin{itemize}
\item $3\beta_2\ge\min\{3,2+\alpha\}$
\item $\alpha\ge\min\{1,2\beta_2-1\}$.
\end{itemize}
We will now prove (\ref{eq:beta21}) and (\ref{eq:alpha}) by contradiction. Assume that at least one of (\ref{eq:beta21}) and
(\ref{eq:alpha}) is false. Then both (\ref{eq:beta21}) and (\ref{eq:alpha}) are false according to the above inequalities.
We see that
\begin{itemize}
\item $3\beta_2\ge2+\alpha$
\item $\alpha\ge2\beta_2-1$
\end{itemize}
which implies that $\frac23+\frac13\alpha\le\frac12+\frac12\alpha$, hence $\alpha\ge1$, a
contradiction. This completes the proof of the relations (\ref{eq:gamma21})--(\ref{eq:alpha}).

After a suitable automorphism of $L[[s]]$, we may assume that
$$
\mathbf a_1(s)=s^{\alpha_1}.
$$
Generically, each arc lives in $N_2$, and thus satisfies the equations of $N_2(8)$. Let $\tilde c$ denote the coefficient of $s^{\alpha_1}$ in the formal power series $\mathbf c_4(s)-i\mathbf a_2(s)^2$ (a priori, $\tilde c$ may or may not be zero). Expanding the equations $\bar f_{2,6}(\mathbf a(s),\mathbf b(s),\mathbf c(s))$, $\bar f_{2,7}(\mathbf a(s),\mathbf b(s),\mathbf c(s))$, $\bar f_{2,8}(\mathbf a(s),\mathbf b(s),\mathbf c(s))$ as power series in $s$ gives:
 \begin{eqnarray*}
0=\bar f_{2,6}(\mathbf a(s),\mathbf b(s),\mathbf c(s))=g_{2,6}s^{3\alpha_1}+h_{2,6}s^{3\alpha_1+1}\\
0=\bar f_{2,7}(\mathbf a(s),\mathbf b(s),\mathbf c(s))=g_{2,7}s^{2\alpha_1}+h_{2,7}s^{2\alpha_1+1}\\
0=\bar f_{2,8}(\mathbf a(s),\mathbf b(s),\mathbf c(s))=g_{2,8}s^{\alpha_1}+h_{2,8}s^{\alpha_1+1},
\end{eqnarray*}
where $g_{2,6},g_{2,7},g_{2,8}$ are polynomials in $a_{np}$, $b_{np}$, $c_{np}$ and $h_{2,6},h_{2,7},h_{2,8}\in L[[s]]$.

Since $\bar f_{2,6}(\mathbf a(s),\mathbf b(s),\mathbf c(s))$, $\bar f_{2,7}(\mathbf a(s),\mathbf b(s),\mathbf c(s))$, $\bar f_{2,8}(\mathbf a(s),\mathbf b(s),\mathbf c(s))$ must vanish identically as power series in $s$, we must have $g_{2,6}=g_{2,7}=g_{2,8}=0$. Let us look at the $g_{2,i}$'s. They are:
\begin{eqnarray}
\ g_{2,6}=\tilde c-2i{a_3}+b_{2,\alpha_1}^3=0\label{eq:g26}\\
\ g_{2,7}=2i({c_5}-6i{a_2}{a_3})+4i\tilde c+3b_{2,\alpha_1}^2{b_3}=0\label{eq:g27}\\
\ g_{2,8}=12i{a_2}^2{a_3}+4i{a_2}{c_5}+2i{a_2}^2\tilde c+3b_{2,\alpha_1}{b_3}^2=0.\label{eq:g28}
\end{eqnarray}
Elements $a_2$, $a_3$, $b_3$, $c_5$ lie in $K(N_4(8))\subset L$ and
are different from 0. Let us regard (\ref{eq:g26})--(\ref{eq:g28}) as
a system of three equations over $L$ in two unknowns
$b_{2,\alpha_1},\tilde c$; if the wedge exists, these equations should
have a solution. Let us prove that this is in fact not the case, thus
obtaining the desired contradiction.

The subfield of $K(N_4(8))$ generated by $a_2$, $a_3$, $b_3$, $c_5$ is
isomorphic to the field of fractions of the ring
$B=\frac{\Bbbk[a_2,a_3,b_3,c_5]}{(b_3^3+2ia_2^2c_5+4a_2^3a_3)}$. Let
$Y$ denote the affine subscheme of $\mathbb A^2_B$ defined by the
equations (\ref{eq:g26})--(\ref{eq:g28}) and let $\bar Y$ denote its
closure in $\mathbb P^2_B$. The scheme $\bar Y$ is defined in $\mathbb
P^2_B$ by the system of three equations
\begin{eqnarray}
\ G_{2,6}=Z^2\tilde C-2i{a_3}Z^3+B_{2,\alpha_1}^3=0\label{eq:g26z}\\
\ G_{2,7}=2i({c_5}-6i{a_2}{a_3})Z^2+4i\tilde
CZ+3B_{2,\alpha_1}^2{b_3}=0\label{eq:g27z}\\
\ G_{2,8}=(12i{a_2}^2{a_3}+4i{a_2}{c_5})Z+2i{a_2}^2\tilde
C+3B_{2,\alpha_1}{b_3}^2=0,\label{eq:g28z}
\end{eqnarray}
homogeneous in the variables $Z,\tilde C,B_{2,\alpha_1}$.

Suppose the system (\ref{eq:g26})--(\ref{eq:g28}) had a solution in
$L$. This means that the natural map $Y\rightarrow\ Spec\ B$ is
dominant, and hence the map $Y\rightarrow\ Spec\ B$ is surjective by
the Proper Mapping Theorem. Thus to prove non-existence of solutions
of (\ref{eq:g26})--(\ref{eq:g28}) it is sufficient to find one
specific $\Bbbk$-rational point of $Spec\ B$ which is not in the image
of $\bar Y$. In other words, it suffices to find specific elements of
$\Bbbk$ such that when these elements are substituted for $a_2$,
$a_3$, $b_3$, $c_5$, the resulting system of homogeneous equations in
$Z,\tilde C,B_{2,\alpha_1}$ has no non-zero solutions. We can easily
find such elements. For example, put
\begin{equation}
b_3=0.\label{eq:b3=0}
\end{equation}
Then
\begin{equation}
2ia_2^2c_5+4a_2^3a_3=0.\label{eq:b_3=0bis}
\end{equation}
We will take
\begin{equation}
a_2\ne0.\label{eq:a2ne0}
\end{equation}
Then equation (\ref{eq:b_3=0bis}) implies that
\begin{equation}
c_5-2ia_2a_3=0.\label{eq:c5=2ia2a3}
\end{equation}
Substituting (\ref{eq:b3=0}) and (\ref{eq:c5=2ia2a3}) into $G_{2,7}$ and $G_{2,8}$, we obtain
\begin{eqnarray}
\bar G_{2,7}=8{a_2}{a_3}Z^2+4i\tilde CZ=0\label{eq:g27zbar}\\
\bar G_{2,8}=(12i-8){a_2}^2{a_3}Z+2i{a_2}^2\tilde C=0.\label{eq:g28zbar}
\end{eqnarray}
\bigskip
If $Z=0$ then, in view of (\ref{eq:a2ne0}) and the equation $G_{2,6}=0$, we have $\tilde C=B_{2,\alpha_1}=0$. Thus there are no non-trivial solutions with $Z=0$. Assume $Z\ne0$ and divide $\bar G_{2,7}$ by $Z$. Now it is easy to see that there exist $a_2,a_3\in\Bbbk$ with $a_2\ne0$ such that the system
\begin{eqnarray}
8{a_2}{a_3}Z+4i\tilde C=0\label{eq:g27zbarnew}\\
\bar G_{2,8}=(12i-8){a_2}^2{a_3}Z+2i{a_2}^2\tilde C=0.\label{eq:g28zbarnew}
\end{eqnarray}
Has
\begin{equation}
Z=\tilde C=0\label{eq:trivial}
\end{equation}
as the only solution. (\ref{eq:trivial}) together with $G_{2,6}$ implies that $B_{2,\alpha_1}=0$. We have proved that there exists a choice of elements $a_2,a_3,b_3,c_5\in\Bbbk$, satisfying $b_3^3+2ia_2^2c_5+4a_2^3a_3=0$, such that after substituting these values into $G_{2,6}=G_{2,7}=G_{2,8}=0$ the resulting system has non non-trivial solutions. This completes the proof of the non-inclusion $\overline{N_4} \not \subset  \overline{N_2}$.

\item  $\bullet${\bf $\overline{N_6} \not \subset \overline{N_4}$.}

In this case we truncate at the order $o_6=12$. We argue by contradiction. Assume that $\overline{N_6(12)}\subset \overline{N_4(12)}$. Let $\phi_{64}$ be a wedge  with generic arc living in $N_4(12)$ and special arc mapping to the generic point of $N_6(12)$. The following equations hold on $N_4(12)$:

\begin{eqnarray*}
\ \mathbf a_1=\mathbf b_1=\mathbf c_1=\mathbf c_2=\mathbf c_3=\mathbf b_2=0\\
\ f_{4,8}=\mathbf c^2_4+\mathbf a_2^4=0\\
\ f_{4,9}=2\mathbf c_4\mathbf c_5+\mathbf b_3^3+4\mathbf a_2^3\mathbf a_3=0\\
\  f_{4,10}=\mathbf c_5^2+2\mathbf c_4\mathbf c_6+3\mathbf b_3^2\mathbf b_4+6\mathbf a_2^2\mathbf a_3^2+4\mathbf a_2^3\mathbf a_4=0\\
\  f_{4,11}=2\mathbf c_5\mathbf c_6+2\mathbf c_4\mathbf c_7+3\mathbf b_3^2\mathbf b_5+3\mathbf b_3\mathbf b_4^2+12\mathbf a_2^2\mathbf a_3\mathbf a_4+4\mathbf a_2^3\mathbf a_5+4\mathbf a_2\mathbf a_3^3=0
\end{eqnarray*}
We have $f_{4,8}=(\mathbf c_4+i\mathbf a_2^2)(\mathbf c_4-i\mathbf a_2^2)$. As can be seen from Figure 1, $\mathcal P_{4,12}$ is the distinguished ideal corresponding to the irreducible factor $\mathbf c_4-i\mathbf a_2^2$. Let us use the notation
$$
g_{2,4}:=\mathbf c_4-i\mathbf a_2^2;
$$
we have $g_{2,4}=0$ on $N_4(12)$.
We have
\begin{equation}
a_{20}=b_{30}=c_{40}=c_{50}=0;\label{eq:vanishonE2}
\end{equation}
we want to show that
\begin{eqnarray}
\mathbf a_2(s)\ &|&\ \mathbf b_3(s)\label{eq:beta3}\\
\mathbf a_2^2(s)\ &|&\ \mathbf c_4(s)\label{eq:gamma4}\\
\mathbf a_2(s)\ &|&\ \mathbf c_5(s)\label{eq:gamma5}
\end{eqnarray}
in $L[[s]]$. \\
The equation $g_{2,4}=0$ implies (\ref{eq:gamma4}).\\
Now, (\ref{eq:beta3}) and (\ref{eq:gamma5}) are equivalent to saying that
\begin{eqnarray}
\alpha_2&\le&\beta_3\quad\text{and}\label{eq:beta3val}\\
\alpha_2&\le&\gamma_5\label{eq:gamma5val}.
\end{eqnarray}
By (\ref{eq:vanishonE2}), we have $\alpha_2>0$.
Using (\ref{eq:gamma4}), equations $f_{4,9}=0$ and
$f_{4,10}=0$ yield
\begin{itemize}
 \item $\beta_3\ge\min\left\{\frac2{3}\alpha_2+\frac1{3}\gamma_5, \alpha_2\right\}$
 \item $\gamma_5\ge\min\{\alpha_2, \beta_3\}$.
\end{itemize}
We prove (\ref{eq:beta3val}) and (\ref{eq:gamma5val}) by contradiction. Suppose at least one of (\ref{eq:beta3val}) and
(\ref{eq:gamma5val}) is false. Then both (\ref{eq:beta3val}) and (\ref{eq:gamma5val}) are false by the above inequalities.
Then
\begin{itemize}
 \item $\beta_3\ge\frac2{3}\alpha_2+\frac1{3}\gamma_5$
 \item $\gamma_5\ge\beta_3$.
\end{itemize}
Hence $\frac23\gamma_5\ge\frac23\alpha_2$, so $\gamma_5\ge\alpha_2$, a contradiction. This completes the proof of
(\ref{eq:beta3})--(\ref{eq:gamma5}).

For the purposes of this non-inclusion, we will deviate slightly from
our standard notation. Namely, we will write $b_3=b_{3\alpha_2}$ and
$c_5=c_{5\alpha_2}$. The meaning of all the other symbols remains
unchanged.

Then the first coefficients of the wedge have to satisfy:

\begin{eqnarray*}
\ c_4-ia_2^2=0\\
\ 2c_4c_5+{b_{3}}^3+4a_2^3{a_3}=0\\
\ {c_5}^2+2c_4{c_6}+3{b_{3}}^2{b_4}+6a_2^2{a_3}^2=0\\
\ 2c_5{c_6}+3b_{3}{b_4}^2+4a_2{a_3}^3=0
\end{eqnarray*}
as well as
\begin{equation}
{c_6}^2+{b_4}^3+{a_3}^4=0.\label{eq:E6forN6}
\end{equation}
Substituting $c_4$ for $ia_2^2$, the above system rewrites as
\begin{eqnarray*}
\ 2ia_2^2c_5+{b_{3}}^3+4a_2^3{a_3}=0\\
\ {c_5}^2+(2i{c_6}+6{a_3}^2)a_2^2+3{b_{3}}^2{b_4}=0\\
\ 2c_5{c_6}+3b_{3}{b_4}^2+4a_2{a_3}^3=0.
\end{eqnarray*}
We view this system as a system of three homogeneous equations over
$L$ in three unknowns $a_2$, $b_3$, $c_5$. The coefficients of the
system are polynomials in $a_3$, $b_4$, $c_6$, which are viewed as
fixed elements of $K(N_6(12))$. Moreover, we must have $a_2\ne0$ by
definition of $a_2$. As in the previous non-inclusion, to prove that
this system has no non-zero solutions, it suffices to find specific
values of $a_3$, $b_4$, $c_6$ in $\Bbbk$ satisfying
(\ref{eq:E6forN6}), such that the resulting system of three equations
has no non-zero solutions. We take $a_3=0$. Then
\begin{equation}
{c_6}^2+{b_4}^3=0\label{eq:a3=0bis}
\end{equation}
and our system becomes
\begin{eqnarray}
\ 2ia_2^2c_5+{b_{3}}^3=0\label{eq:49}\\
\ {c_5}^2+2i{c_6}a_2^2+3{b_{3}}^2{b_4}=0\label{eq:410}\\
\ 2c_5{c_6}+3b_{3}{b_4}^2=0.\label{eq:411}
\end{eqnarray}
We work in a finite extension of $K(N_6(12))$ which contains a square
root of $b_4$; we pick and fix one of the two possible square roots
and denote it by $b_4^{\frac12}$. From (\ref{eq:a3=0bis}) we obtain
\begin{equation}
c_6=-b_4^{\frac32}.\label{eq:c6b43halves}
\end{equation}
Substituting (\ref{eq:c6b43halves}) into (\ref{eq:411}) and dividing
through by $b_4^{\frac32}$, we obtain
\begin{equation}
c_5=\frac32b_3b_4^{\frac12}.\label{eq:c_5=}
\end{equation}
Substituting (\ref{eq:c_5=}) into (\ref{eq:49}) yields
\begin{equation}
b_3^2=-2ib_4^{\frac12}a_2^2.\label{eq:b_32=}
\end{equation}
Finally, substituting (\ref{eq:c_5=}) and (\ref{eq:b_32=}) into
(\ref{eq:410}), we obtain
\begin{equation}
\left(-\frac{27}4-2-9\right)ib_4^{\frac32}a_2^2=0.\label{eq:final}
\end{equation}
Now, substitute suitable non-zero elements of $\Bbbk$ for
$b_4^{\frac12}$ and $c_6$ in such a way that (\ref{eq:a3=0bis}) is
satisfied. By (\ref{eq:final}), any solution of the resulting system
of equations satisfies $a_2=0$. Then $b_3=c_5=0$ from
(\ref{eq:49})--(\ref{eq:411}). Thus our system of equations has no
non-zero solutions, as desired. This completes the proof of the
non-inclusion $\overline{N_6} \not \subset \overline{N_4}$.

\item  $\bullet${\bf $\overline{N_6} \not \subset \overline{N_1}$.}

In this case we truncate at the order $o_6=12$. We argue by
contradiction: suppose that $\overline{N_6}(12)\subset
\overline{N_1}(12)$.
Let $\phi_{61}$ be a wedge with the generic arc living in $N_1$ and
the special arc mapping to the generic point of $N_6(12)$.

The following equations hold on $N_1(12)$:
\begin{eqnarray*}
&\mathbf a_1=\mathbf b_1=\mathbf c_1=\mathbf c_2=0\\
&f_{1,6}=\mathbf c_3^2+\mathbf b_2^3=0\\
&f_{1,7}= 2\mathbf c_3\mathbf c_4+3\mathbf b_2^2\mathbf b_3=0\\
&f_{1,8}=\mathbf c_4^2+2\mathbf c_3\mathbf c_5+3\mathbf b_2^2\mathbf
b_4+3\mathbf b_2\mathbf b_3^2+\mathbf a_2^4=0\\
&f_{1,9}=2\mathbf c_4\mathbf c_5+2\mathbf c_3\mathbf c_6+\mathbf
b_3^3+6\mathbf b_2\mathbf b_3\mathbf b_4+3\mathbf b_2^2\mathbf
b_5+4\mathbf a_2^3\mathbf a_3=0\\
&f_{1,10}=\mathbf c_5^2+2\mathbf c_4\mathbf c_6+2\mathbf c_3\mathbf
c_7+3\mathbf b_3^2\mathbf b_4+3\mathbf b_2^2\mathbf b_6+6\mathbf
b_2\mathbf b_3\mathbf b_5+3\mathbf b_2\mathbf b_4^2+6\mathbf
a_2^2\mathbf a_3^2+4\mathbf a_2^3\mathbf a_4=0\\
&f_{1,11}=2\mathbf c_5\mathbf c_6+2\mathbf c_4\mathbf c_7+2\mathbf
c_3\mathbf c_8+3\mathbf b_2^2\mathbf b_7+3\mathbf b_3^2\mathbf
b_5+3\mathbf b_3\mathbf b_4^2+6\mathbf b_2\mathbf b_3\mathbf
b_6+6\mathbf b_2\mathbf b_4\mathbf b_5+4\mathbf a_2^3\mathbf a_5+\\
&+12\mathbf a_2^2\mathbf a_3\mathbf a_4+4\mathbf a_2\mathbf a_3^3=0
\end{eqnarray*}
The following equations come from the equations of $N_6(12)$:
$$
a_{10}=a_{20}=b_{10}=b_{20}=b_{30}=c_{10}=c_{20}=c_{30}=c_{40}=c_{50}=0.
$$
We want to prove the following
divisibility relations:
\begin{eqnarray}
\mathbf b_2(s)\ &|&\ \mathbf a_2(s)^2\label{eq:alpha2}\\
\mathbf b_2(s)\ &|&\ \mathbf b_3(s)^2\label{eq:beta31}\\
\mathbf b_2(s)^3\ &|&\ \mathbf c_3(s)^2\label{eq:gamma32}\\
\mathbf b_2(s)\ &|&\ \mathbf c_4(s)\label{eq:gamma41}\\
\mathbf b_2(s)\ &|&\ \mathbf c_5(s)^2.\label{eq:gamma51}
\end{eqnarray}
To do this, it is sufficient to show that
\begin{eqnarray}
\gamma_3&=&\frac{3}{2}\beta_2\label{eq:gamma32val}\\
\gamma_4&\ge&\beta_2\label{eq:gamma41val}\\
\alpha_2,\beta_3,\gamma_5&\ge&\frac{1}{2}\beta_2.\label{eq:alpha2val}
\end{eqnarray}
We have $\beta_2>0$.
The equality (\ref{eq:gamma32val}) is immediate from
$f_{1,6}=0$. (\ref{eq:gamma41val}) follows from
$f_{1,6}=f_{1,7}=0$ and (\ref{eq:alpha2val}). It remains to prove (\ref{eq:alpha2val}), which is equivalent to saying
that
\begin{equation}
\min\{\alpha_2,\beta_3,\gamma_5\}\ge\frac{1}{2}\beta_2.\label{eq:min1half}
\end{equation}
We prove (\ref{eq:min1half}) by contradiction. Let $M=\min\{\alpha_2,\beta_3,\gamma_5\}$ and assume that
\begin{equation}
M<\frac12\beta_2.\label{eq:lessthanhalf}
\end{equation}
Equations $f_{1,6}=f_{1,7}=0$ can be interpreted as saying that $\frac{\mathbf c_3(s)}{\mathbf b_2(s)^{\frac32}}$ and
$\frac{\mathbf c_4(s)}{\mathbf b_2(s)^{\frac12}\mathbf b_3(s)}$ is invertible in a suitable finite extension $B$ of $L[[s]]$. Substituting $\mathbf c_3(s)$ and $\mathbf c_4(s)$ in $f_{1,8}$,
$f_{1,9}$ and $f_{1,10}$ by suitable multiples of $\mathbf b_2(s)^{\frac32}$ and
$\mathbf b_2(s)^{\frac12}\mathbf b_3(s)$ by a unit of $B$,
we obtain the following inequalities:
\begin{eqnarray}
\alpha_2&\ge&\frac{1}{4}\min\left\{\frac32\beta_2+\gamma_5,2\beta_2,\beta_2+2\beta_3\right\}\label{eq:alpha2>}\\
\beta_3&\ge&\frac{1}{3}\min\left\{\frac12\beta_2+\beta_3+\gamma_5,\beta_2+\beta_3,\frac32\beta_2,3\alpha_2\right\}\label{eq:beta3>}\\
\gamma_5&\ge&\frac{1}{2}\min\left\{\frac12\beta_2+\beta_3,2\beta_3,\beta_2,2\alpha_2\right\}.\label{eq:gamma5>}
\end{eqnarray}
Now, (\ref{eq:lessthanhalf}), (\ref{eq:alpha2>}) and the definition of $M$ imply that
\begin{equation}
M<\alpha_2\label{eq:M<a2}
\end{equation}
(indeed, if we had $M\ge\alpha_2$, we could use (\ref{eq:lessthanhalf}) and the definition of $M$ to show that $M$ is strictly less than each of the three quantities on the right hand side of (\ref{eq:alpha2>}), which would be a contradiction).

In a similar way, (\ref{eq:lessthanhalf}), (\ref{eq:beta3>}), (\ref{eq:M<a2}) and the definition of $M$ imply that
\begin{equation}
M<\beta_3.\label{eq:M<b3}
\end{equation}
By (\ref{eq:M<a2}) and (\ref{eq:M<b3}), we have $M=\gamma_5$, which contradicts (\ref{eq:gamma5>}) (using (\ref{eq:M<a2})
and (\ref{eq:M<b3}) once again). This completes the proof of (\ref{eq:alpha2})--(\ref{eq:gamma51}).

\bigskip


Replacing $s$ by $s^2$ in the parametrization of the wedge, we may
assume, without loss of generality, that $\beta_2$ is even. The first
coefficients of the wedge must satisfy the following equations (as
above we change the notation by $c_4=c_{4,\beta_2}$,
$c_5=c_{5,\frac{\beta_2}2}$, $b_3=b_{3,\frac{\beta_2}2}$ and
$a_2=a_{2,\frac{\beta_2}2}$):
\begin{eqnarray*}
&\  c_3^2+\  b_2^3=0\\
& 2\  c_3\  c_4+3\  b_2^2\  b_3=0\\
&\  c_4^2+2\  c_3\  c_5+3\  b_2^2\  b_4+3\  b_2\  b_3^2+\  a_2^4=0\\
&2\  c_4\  c_5+2\  c_3\  c_6+\  b_3^3+6\  b_2\  b_3\  b_4+4\  a_2^3\  a_3=0\\
&\  c_5^2+2\  c_4\  c_6+3\  b_3^2\  b_4+3\  b_2\  b_4^2+6\  a_2^2\  a_3^2=0\\
&2\  c_5\  c_6+3\  b_3\  b_4^2+4\  a_2\  a_3^3=0
\end{eqnarray*}

as well as
\begin{equation}
{c_6}^2+{b_4}^3+{a_3}^4=0.\label{eq:E6forN6bis}
\end{equation}

We view this system as a system of six homogeneous equations over
$L$ in six unknowns $a_2$, $b_2$, $b_3$, $c_3$, $c_4$,  $c_5$. The
coefficients of the system are polynomials in $a_3$, $b_4$, $c_6$,
which are viewed as fixed elements of $K(N_6(12))$. As in the previous
non-inclusion, to prove that this system has no non-zero solutions, it
suffices to find specific values of $a_3$, $b_4$, $c_6$ in $\Bbbk$
satisfying (\ref{eq:E6forN6}), such that the resulting system of six
equations has no non-zero solutions. In this case, we take $a_3=0$,
$c_6=1$ and $b_4$ a non-real root of $z^3=-1$. We obtain:
\begin{eqnarray*}
&\  c_3^2+\  b_2^3=0\\
& 2\  c_3\  c_4+3\  b_2^2\  b_3=0\\
&\  c_4^2+2\  c_3\  c_5+3(1/2-\sqrt3/2i)\  b_2^2\  +3\  b_2\
b_3^2+\  a_2^4=0\\
&2\  c_4\  c_5+2\  c_3+\  b_3^3+6(1/2-\sqrt3/2i)\  b_2\  b_3\  =0\\
&\  c_5^2+2\  c_4+3(1/2-\sqrt3/2i)\  b_3^2+3(1/2-\sqrt3/2i)^2\  b_2\ =0\\
&2\  c_5+3(1/2-\sqrt3/2i)^2\  b_3=0
\end{eqnarray*}

Then we ask Maple to solve it and the solution that Maple gives is:
$\{c5 = 0, b3 = 0, a2 = 0, c4 = 0, c3 = 0, b2 = 0\}$,
so that the unique solution is the zero one.

\item  $\bullet${\bf $\overline{N_6} \not \subset \overline{N_2}$.}

In this case we truncate at the order $o_6=12$. We argue by
contradiction: suppose that $\overline{N_6}(12)\subset
\overline{N_2}(12)$.
Let $\phi_{62}$ be a wedge with the generic arc living in $N_2$ and
the special arc mapping to the generic point of $N_6(12)$. The following
equations vanish on $N_2(12)$:
\begin{eqnarray*}
\mathbf b_1=\mathbf c_1=0\\
\ g_{2,2}= \mathbf c_2-i\mathbf a_1^2=0 \\
\  \bar f_{2,3}=\mathbf c_3-2i\mathbf a_1\mathbf a_2=0\\
\ f_{2,6}=\mathbf c_3^2+2\mathbf c_2\mathbf c_4+\mathbf b_2^3+4\mathbf
a_1^3\mathbf a_3+6\mathbf a_1^2\mathbf a_2^2=0\\
\ f_{2,7}= 2\mathbf c_3\mathbf c_4+2\mathbf c_2\mathbf c_5+3\mathbf
b_2^2\mathbf b_3+4\mathbf a_1^3\mathbf a_4+12 \mathbf a_1^2\mathbf
a_2\mathbf a_3+4\mathbf a_2^3\mathbf a_1=0\\
\vdots \\
\  f_{2,11}=3\mathbf b_2^2\mathbf b_7+...+2\mathbf c_2\mathbf
c_9+...+4\mathbf a_2^3\mathbf a_5+...=0;
\end{eqnarray*}
We write it in the following way:
\begin{eqnarray*}
\ &\mathbf b_1=\mathbf c_1=0\\
\ &g_{2,2}= \mathbf c_2-i\mathbf a_1^2=0 \\
\  &\bar f_{2,3}=\mathbf c_3-2i\mathbf a_1\mathbf a_2=0\\
\  &f_{2,6}=\mathbf b_2^3+2i\mathbf a_1^2(\mathbf c_4-i\mathbf
a_2^2-2i\mathbf a_1\mathbf a_3)=0\\
\ &f_{2,7}=3\mathbf b_2^2\mathbf b_3+2i\mathbf a_1^2(\mathbf
c_5-2i\mathbf a_2\mathbf a_3-2i\mathbf a_1\mathbf a_4)+4i\mathbf
a_1\mathbf a_2(\mathbf c_4-i\mathbf a_2^2-2i\mathbf a_1\mathbf
a_3)=0\\
&f_{2,8}=3\mathbf b_2^2\mathbf b_4+3\mathbf b_2\mathbf b_3^2+2i\mathbf
a_1^2(\mathbf c_6-i\mathbf a_3^2-2i\mathbf a_2\mathbf a_4-2i\mathbf
a_1\mathbf a_5)+\\
&+4i\mathbf a_1\mathbf a_2(\mathbf c_5-2i\mathbf a_2\mathbf
a_3-2i\mathbf a_1\mathbf a_4)+(\mathbf c_4+i\mathbf a_2^2+2i\mathbf
a_1\mathbf a_3)(\mathbf c_4-i\mathbf a_2^2-2i\mathbf a_1\mathbf
a_3)=0\\
&f_{2,9}=\mathbf b_3^3+3\mathbf b_2^2\mathbf b_5+6\mathbf b_2\mathbf
b_3\mathbf b_4+2i\mathbf a_1^2(\mathbf c_7-2i\mathbf a_3\mathbf
a_4-2i\mathbf a_2\mathbf a_5-2i\mathbf a_1\mathbf a_6)+\\
&+4i\mathbf a_1\mathbf a_2(\mathbf c_6-i\mathbf a_3^2-2i\mathbf
a_2\mathbf a_4-2i\mathbf a_1\mathbf a_5)+(\mathbf c_4+i\mathbf
a_2^2+2i \mathbf a_1\mathbf a_3)(\mathbf c_5-2i\mathbf a_2\mathbf
a_3-2i\mathbf a_1\mathbf a_4)+\\
&+(\mathbf c_5+2i\mathbf a_2\mathbf a_3+2i\mathbf a_1\mathbf
a_4)(\mathbf c_4-i\mathbf a_2^2-2i\mathbf a_1\mathbf a_3)=0\\
&f_{2,10}=3\mathbf b_2^2\mathbf b_6+3\mathbf b_3^2\mathbf b_4+3\mathbf
b_4^2\mathbf b_2+6\mathbf b_2\mathbf b_3\mathbf b_5+2i\mathbf
a_1^2(\mathbf c_8-i\mathbf a_4^2-2i\mathbf a_3\mathbf a_5-2i\mathbf
a_2\mathbf a_6-2i\mathbf a_1\mathbf a_7)\\
&+4i\mathbf a_1\mathbf a_2(\mathbf c_7-2i\mathbf a_3\mathbf
a_4-2i\mathbf a_2\mathbf a_5-2i\mathbf a_1\mathbf a_6)+\\
&+(\mathbf c_4+i\mathbf a_2^2+2i\mathbf a_1\mathbf a_3)(\mathbf
c_6-i\mathbf a_3^2-2i\mathbf a_2\mathbf a_4-2i\mathbf a_1\mathbf
a_5)+\\
&+(\mathbf c_5+2i\mathbf a_2\mathbf a_3+2i\mathbf a_1\mathbf
a_4)(\mathbf c_5-2i\mathbf a_2\mathbf a_3-2i\mathbf a_1\mathbf a_4)+\\
&+(\mathbf c_6+i\mathbf a_3^2+2i\mathbf a_2\mathbf a_4+2i\mathbf
a_1\mathbf a_5)(\mathbf c_4-i\mathbf a_2^2-2i\mathbf a_1\mathbf
a_3)=0\\
& f_{2,11}=3\mathbf b_2^2\mathbf b_7+6\mathbf b_2\mathbf b_3\mathbf
b_6+6\mathbf b_2\mathbf b_4\mathbf b_5+3\mathbf b_3\mathbf b_4^2+\\
&+3\mathbf b_3^2\mathbf b_5+2i\mathbf a_1^2(\mathbf c_9-2i\mathbf
a_4\mathbf a_5-2i\mathbf a_3\mathbf a_6-2i\mathbf a_2\mathbf
a_7-2i\mathbf a_1\mathbf a_8)\\
&+4i\mathbf a_1\mathbf a_2(\mathbf c_8-i\mathbf a_4^2-2i\mathbf
a_3\mathbf a_5-2i\mathbf a_2\mathbf a_6-2i\mathbf a_1\mathbf a_7)+\\
&+(\mathbf c_4+i\mathbf a_2^2+2i\mathbf a_1\mathbf a_3)(\mathbf
c_7-2i\mathbf a_3\mathbf a_4-2i\mathbf a_2\mathbf a_5-2i\mathbf
a_1\mathbf a_6)\\
&+(\mathbf c_5+2i\mathbf a_2\mathbf a_3+2i\mathbf a_1\mathbf
a_4)(\mathbf c_6-i\mathbf a_3^2-2i\mathbf a_2\mathbf a_4-2i\mathbf
a_1\mathbf a_5)\\
&+(\mathbf c_6+i\mathbf a_3^2+2i\mathbf a_2\mathbf a_4+2i\mathbf
a_1\mathbf a_5)(\mathbf c_5-2i\mathbf a_2\mathbf a_3-2i\mathbf
a_1\mathbf a_4)\\
&+(\mathbf c_7+2i\mathbf a_3\mathbf a_4+2i\mathbf a_2\mathbf
a_5+2i\mathbf a_1\mathbf a_6)(\mathbf c_4-i\mathbf a_2^2-2i\mathbf
a_1\mathbf a_3)=0.\\
\end{eqnarray*}

The following equations come from the equations of $N_6(12)$:
$$
a_{10}=a_{20}=b_{10}=b_{20}=b_{30}=c_{10}=c_{20}=c_{30}=c_{40}=c_{50}=0.
$$
In this case, because of the number of variables, the computation is
more difficult than for the other cases. We want to compute or at least bound
below the rational numbers

\begin{eqnarray}
\alpha'_2 &:=&\ \frac{\alpha_2}{\alpha_1}\label{eq:alpha2bis}\\
\beta'_2\ &:=&\ \frac{\beta_2}{\alpha_1}\label{eq:beta31bis}\\
\beta'_3\ &:=&\ \frac{\beta_3}{\alpha_1}\label{eq:gamma32bis}\\
\gamma'_4\ &:=&\ \frac{\gamma_4}{\alpha_1}\label{eq:gamma41bis}\\
\gamma'_5\ &:=&\ \frac{\gamma_5)}{\alpha_1}.\label{eq:gamma51bis}
\end{eqnarray}

We use the following dichotomy.\\
If $\alpha_2\ge \frac{1}{2}\alpha_1$ then we have

\begin{eqnarray}
{\beta_2}\ge\alpha_1\\
{\beta_3}\ge \ \frac{1}{2}\alpha_1\\
{\gamma_4}\ \ge\alpha_1\\
{\gamma_5}\ge \frac{1}{2}\alpha_1.
\end{eqnarray}

We try to construct a wedge as usual. Replacing $s$ by $s^2$ in the
parametrization of the wedge, we may assume, without loss of
generality, that $\alpha_1$ is even. We deviate from our standard
notation (only for the purposes of the case $\alpha_2\ge \frac{1}{2}\alpha_1$), in that we put $a_2=a_{2,\frac{\alpha_1}2}$
$b_2=b_{2,\alpha_1}$, $b_3=b_{3,\frac{\alpha_1}2}$,
$c_4=c_{4,\alpha_1}$, $c_5=c_{5,\frac{\alpha_1}2}$.

In this case the first equations are of the form
\begin{eqnarray*}
\  & b_2^3+2i a_1^2( c_4-ia_2^2-2i a_1 a_3)=0\\
\ &3 b_2^2 b_3+2i a_1^2(c_5-2i a_2 a_3)+4i a_1 a_2( c_4-i a_2^2-2i
a_1a_3)=0\\
\ &3 b_2^2 b_4+3 b_2 b_3^2+2ia_1^2( c_6-i a_3^2)+4i a_1 a_2( c_5-2i
a_2a_3)+\\
&+( c_4+i a_2^2+2ia_1 a_3)( c_4-i a_2^2-2i a_1a_3)=0\\
\ & b_3^3+6 b_2b_3 b_4+4i a_1 a_2( c_6-i a_3^2)+( c_4+ia_2^2+2i a_1
a_3)( c_5-2i a_2a_3)+\\
&+( c_5+2i a_2 a_3)( c_4-i a_2^2-2i a_1a_3)=0\\
\ &3 b_4^2 b_2+3 b_3^2 b_4+( c_4+i a_2^2+2i a_1 a_3)( c_6-i a_3^2)+\\
&+(c_5+2i a_2 a_3)( c_5-2i a_2 a_3)+( c_6+i a_3^2)( c_4-i a_2^2-2i a_1
a_3)=0\\
\ & 3 b_3 b_4^2+( c_5+2i a_2 a_3)( c_6-i a_3^2)+( c_6+i a_3^2)( c_5-2i
a_2 a_3)=0.
\end{eqnarray*}
Thanks to XMaple, taking in this case $b_4=0$ and $a_3=1$ one can show
that the above system of equations, combined with the first equation
of $N_6$,
$$
a_3^4+b_4^3+c_6^2=0,
$$
has no non-zero solutions, so the wedge cannot be constructed.
\bigskip

From now on we shall assume that $\alpha_2 < \frac{1}{2}\alpha_1$. One always has $\beta_2
\ge \frac{2}{3}\alpha_1$ thanks to the equation $f_{2,6}$.\smallskip

For each equation, let us write the $\mu$-adic orders of monomials
appearing in it, which can possibly be the lowest for this equation:\smallskip

\begin{eqnarray*}
\ f_{2,6}:\ \ 3\beta_2,\ 2\alpha_1+\gamma_4,\ 2\alpha_1+2\alpha_2\\
\ f_{2,7}:\ 2\beta_2+\beta_3,\ 2\alpha_1+\gamma_5,\
2\alpha_1+\alpha_2,\ \alpha_1+\alpha_2+\gamma_4,\ \alpha_1+3\alpha_2\\
\ f_{2,8}:\ 2\beta_2,\ \beta_2+2\beta_3,\ 2\alpha_1,\
\alpha_1+\alpha_2+\gamma_5,\ \alpha_1+2\alpha_2,\ 2\gamma_4,\ 4\alpha_2\\
\ f_{2,9}:\ 3\beta_3,\ 2\beta_2,\  \beta_2+\beta_3,\
\gamma_5+\gamma_4,\ 2\alpha_2+\gamma_5,\ \gamma_4+\alpha_2,\ 3\alpha_2\\
\ f_{2,10}:\ \beta_2,\ 2\beta_3,\ \gamma_4,\
2\alpha_2,\ 2\gamma_5\\
\ f_{2,11}:\ \beta_2,\ \beta_3,\ \gamma_4,\ \alpha_2,\
\gamma_5\\
\end{eqnarray*}
\noi\textbf{Note:} Here we have used the fact (easy to prove) that the
following four expressions have $\mu$-adic value equal to zero:
$\mathbf c_6-i\mathbf a_3^2-2i\mathbf a_2\mathbf a_4-2i\mathbf
a_1\mathbf a_5$, $\mathbf c_6+i\mathbf a_3^2+2i\mathbf a_2\mathbf
a_4+2i\mathbf a_1\mathbf a_5$, $\mathbf c_7-2i\mathbf a_3\mathbf
a_4-2i\mathbf a_2\mathbf a_5-2i\mathbf a_1\mathbf a_6$, $\mathbf
c_7+2i\mathbf a_3\mathbf a_4+2i\mathbf a_2\mathbf a_5+2i\mathbf
a_1\mathbf a_6$.
\medskip

Suppose that $\beta_3 \ge \frac{1}{2}\alpha_1(>\alpha_2)$. \smallskip

Now, if $\gamma_4\le \alpha_2$ then from the equation $f_{2,11}$ we
see that $\gamma_5\ge\gamma_4$ (otherwise the term with
$\mu$-adic value $\gamma_5$ would be the only dominant term). But then
the term of value $\gamma_4$ is the only dominant term in $f_{2,10}$, a
contradiction; so
$$
\gamma_4 > \alpha_2.
$$
\medskip

Then the $\mu$-adic values of possible dominant terms are:\smallskip

\begin{eqnarray*}
\ f_{2,6}:\ \ 3\beta_2,\ 2\alpha_1+\gamma_4,\ 2\alpha_1+2\alpha_2\\
\ f_{2,7}:\ 2\beta_2+\beta_3,\  2\alpha_1+\gamma_5,\
\alpha_1+\alpha_2+\gamma_4,\ \alpha_1+3\alpha_2\\
\ f_{2,8}:\ 2\beta_2,\ \beta_2+2\beta_3,\ \alpha_1+\alpha_2+\gamma_5,\
2\gamma_4,\ 4\alpha_2\\
\ f_{2,9}:\ 3\beta_3,\ 2\beta_2,\ \beta_2+\beta_3,\
3\alpha_2,\ \gamma_5+\gamma_4,\ 2\alpha_2+\gamma_5,\
\gamma_4+\alpha_2\\
\ f_{2,10}:\ \beta_2,\ \gamma_4,\ 2\alpha_2,\ 2\gamma_5\\
\ f_{2,11}:  \alpha_2,\ \gamma_5.\\
\end{eqnarray*}
From $f_{2,11}$ we have $\alpha_2=\gamma_5$.
\begin{itemize}
\item First case :\smallskip

Suppose that $\gamma_4>2\alpha_2$.\smallskip

Then the dominant values are:\smallskip

\begin{eqnarray*}
\ f_{2,6}:\ \ 3\beta_2,\ 2\alpha_1+2\alpha_2\\
\ f_{2,7}:\ 2\beta_2+\beta_3,\ \alpha_1+3\alpha_2\\
\ f_{2,8}:\ 2\beta_2,\ \beta_2+2\beta_3,\ 4\alpha_2\\
\ f_{2,9}:\ 2\beta_2,\  \beta_2+\beta_3,\ 2\alpha_2+\gamma_5=3\alpha_2\\
\ f_{2,10}:\ \beta_2,\ 2\alpha_2,\ 2\gamma_5\\
\ f_{2,11}: \ \alpha_2,\ \gamma_5\\
\end{eqnarray*}

So by $f_{2,6}$ we have : $\ 3\beta_2=2\alpha_1+2\alpha_2$.\smallskip

And by  $f_{2,7}$ we have :
$$
2\beta_2+\beta_3=\alpha_1+3\alpha_2.
$$
The two equations imply that $\beta_3<\frac{1}{2}\alpha_1$, a
contradiction.\medskip

\item second case:\smallskip

Suppose that $\gamma_4<2\alpha_2$.\smallskip

Then the possible dominant values are:\smallskip

\begin{eqnarray*}
\ f_{2,6}:\ \ 3\beta_2,\ 2\alpha_1+\gamma_4\\
\ f_{2,7}:\ 2\beta_2+\beta_3,\ \alpha_1+\gamma_4+\alpha_2\\
\ f_{2,8}:\ 2\beta_2,\ \beta_2+2\beta_3,\ 2\gamma_4\\
\ f_{2,9}:\   2\beta_2,\ \beta_2+\beta_3,\ \alpha_2+\gamma_4\\
\ f_{2,10}:\ \beta_2,\ \gamma_4\\
\ f_{2,11}: \ \alpha_2,\ \gamma_5\\
\end{eqnarray*}

So by $f_{2,6}$ and $f_{2,10}$ we have : $\beta_2=\alpha_1=\gamma_4$, a
contradiction (as $\alpha_2<\frac{1}{2}\alpha_1$).\medskip

\item Last case : $\gamma_4=2\alpha_2$.

The dominant values are:
\begin{eqnarray*}
\ f_{2,6}:\ \ 3\beta_2=\ 2\alpha_1+\gamma_4= 2\alpha_1+2\alpha_2\\
\ f_{2,7}:\ \alpha_2+\gamma_4= \alpha_1+3\alpha_2\\
\ f_{2,8}:\  2\gamma_4=4\alpha_2\\
\ f_{2,9}:\
\gamma_5+\gamma_4=2\alpha_2+\gamma_5=\gamma_4+\alpha_2=3\alpha_2\\
\ f_{2,10}:\  \gamma_4=2\alpha_2=2\gamma_5\\
\ f_{2,11}:\  \alpha_2=\gamma_5\\
\end{eqnarray*}
 Then the first equations of the wedge are:\smallskip

 \begin{eqnarray*}
\  &\bar f_{2,6}= b_2^3+2i a_1^2( c_4-i a_2^2)=0\\
\ &\bar f_{2,7}=4i a_1 a_2( c_4-i a_2^2)=0\\
&\bar f_{2,8}=( c_4+i a_2^2)( c_4-i a_2^2)=0\\
&\bar f_{2,9}=\ ( c_4+i a_2^2)( c_5-2i a_2 a_3)+( c_5+2i a_2 a_3)( c_4-i
a_2^2)=0\\
&\bar f_{2,10}=( c_4+i a_2^2)( c_6-i a_3^2)+( c_5+2i a_2 a_3)( c_5-2i a_2
a_3)+( c_6+i a_3^2)( c_4-i a_2^2)=0\\
&\bar f_{2,11}=( c_5+2i a_2 a_3)( c_6-i a_3^2)+( c_6+i a_3^2)( c_5-2i a_2
a_3)=0.\\
& f_{6,12}=a_3^4+b_4^3+c_6^2=0
\end{eqnarray*}
 \end{itemize}

By definitions, we are looking for solutions with $a_1$, $a_2$, $b_2$
different from 0. It is easy to see that this is impossible already
from the equations $\bar f_{2,6}$ and $\bar f_{2,7}$. Indeed, since
$a_1,a_2\ne0$ we have $c_4-i a_2^2=0$ by $\bar f_{2,7}$. Then $\bar
f_{2,6}$ shows that $b_2=0$. This completes the proof of the
non-existence of the wedge in the case $\beta_3\ge\frac{1}{2}\alpha_1$.
\bigskip

Thus we will assume from now on that $\beta_3 <\frac{1}{2}\alpha_1$. Then
$\beta_2> \beta_3$ and the possible dominant values are:\smallskip

\begin{eqnarray*}
\ f_{2,6}:\ \ 3\beta_2,\ 2\alpha_1+\gamma_4,\ 2\alpha_1+2\alpha_2\\
\ f_{2,7}:\ 2\beta_2+\beta_3,\ 2\alpha_1+\gamma_5,\
\ \alpha_1+\alpha_2+\gamma_4,\ \alpha_1+3\alpha_2\\
\ f_{2,8}:\ 2\beta_2,\ \beta_2+2\beta_3,\
\alpha_1+\alpha_2+\gamma_5,\ 2\gamma_4,\ 4\alpha_2\\
\ f_{2,9}:\ 3\beta_3,\  \beta_2+\beta_3,\
\ \gamma_5+\gamma_4,\ 2\alpha_2+\gamma_5,\
\gamma_4+\alpha_2,\  \ 3\alpha_2\\
\ f_{2,10}:\ \beta_2,\ 2\beta_3,\  \gamma_4,\
2\alpha_2,\ 2\gamma_5\\
\ f_{2,11}:\ \ \beta_3,\ \gamma_4,\ \alpha_2,\
\gamma_5\\
\end{eqnarray*}
\begin{itemize}
\item Suppose that $\gamma_5\ge \frac{1}{2}\alpha_1$.

\begin{itemize}
\item If
\begin{equation}
\gamma_4\le \alpha_2\label{eq:gamma4<alpha2}
\end{equation}
then $\beta_3\le \alpha_2$ by
  $f_{2,11}$. Hence $\gamma_4$ becomes the only dominant value in
  $f_{2,10}$ which is not possible. Thus $\gamma_4>\alpha_2$, which
  implies that $\beta_3=\alpha_2$ by $f_{2,11}$.
\item If
\begin{equation}
\gamma_4< 2 \alpha_2,\label{eq:gamma4<2alpha2}
\end{equation}
then $\beta_2=\gamma_4$ by $f_{2,10}$
  and hence $\beta_2=\alpha_1$ by $f_{2,6}$. Thus $\gamma_4=\alpha_1$,
  which contradicts (\ref{eq:gamma4<2alpha2}) and the fact that
  $\alpha_2<\frac12\alpha_1$.
\item If $\gamma_4 >2\alpha_2$ then $\beta_2=2\alpha_2$ by $f_{2,10}$
  and hence $\alpha_2=\frac25\alpha_1$, $\beta_2=\frac45\alpha_1$ by
  $f_{2,6}$. Using the fact that $\gamma_5>0$, we see that in
  $f_{2,7}$, $2\beta_2+\alpha_2$ is the only dominant value, a contradiction.
\item The remaining case is $\gamma_4=2\alpha_2$. The dominant values are:
\begin{eqnarray*}
\ f_{2,6}:\ \ 3\beta_2\ge \ 2\alpha_1+\gamma_4=\ 2\alpha_1+2\alpha_2\\
\ f_{2,7}:\ 2\beta_2+\beta_3\ge
\ \alpha_1+\alpha_2+\gamma_4= \alpha_1+3\alpha_2\\
\ f_{2,8}:\ 2\beta_2,\ \beta_2+2\beta_3,\
 2\gamma_4=\ 4\alpha_2\\
\ f_{2,9}:\ 3\beta_3,\ \beta_2+\beta_3,\
\gamma_4+\alpha_2= \ 3\alpha_2\\
\ f_{2,10}:\ \beta_2\ge2\beta_3= \gamma_4\
=2\alpha_2\\
\ f_{2,11}:\ \ \beta_3=\ \alpha_2,\
\\
\end{eqnarray*}
If $\beta_2\le 2\alpha_2$ then by $f_{2,6}$ we would have
$\beta_2\ge\alpha_1$ and hence $\alpha_2\ge \frac{1}{2}\alpha_1$, a
contradiction. Thus $\beta_2>2\alpha_2$.
\smallskip

\noi\textbf{Claim.} The only dominant values in $f_{2,7}$ are
$$
\alpha_1+\alpha_2+\gamma_4= \alpha_1+3\alpha_2
$$
\noi\textit{Proof of Claim} If not, we would have
$$
2\beta_2+\beta_3= 2\beta_2+\alpha_2= \ \alpha_1+3\alpha_2
$$
Then $2\beta_2=\alpha_1+2\alpha_2\le 3\beta_2-\alpha_1$
(by $f_{2,6}$), thus $ \beta_2\ge \alpha_1$.

We obtain  $2\beta_2+\beta_3\ge
2\alpha_1+\beta_3=2\alpha_1+\alpha_2>\alpha_1+3\alpha_2 $ a
contradiction. This proves the Claim.
\smallskip

Then the first two equations of the wedge are

 \begin{eqnarray*}
\  &\bar f_{2,6}= b_2^3+2i a_1^2( c_4-i a_2^2)=0\\
\ &\bar f_{2,7}=4i a_1 a_2( c_4-i a_2^2)=0\\
\end{eqnarray*}
so there are no solutions with $b_2\ne0$, $a_2\ne0$, $a_1\ne0$,
contradiction.
\end{itemize}

\item Thus $ \gamma_5<\frac{1}{2}\alpha_1$.

First of all, we claim that $\gamma_4$ cannot be dominant in
$f_{2,11}$. Indeed, suppose it was, in other words, suppose that
$\gamma_4\le\min\{\beta_3,\alpha_2,\gamma_5\}$. In particular,
\begin{equation}
\gamma_4<\frac12\alpha_1.\label{eq:gamma4<onehalf}
\end{equation}
Then by $f_{2,10}$ we have
\begin{equation}
\beta_2=\gamma_4.\label{eq:beta2=gamma_4}
\end{equation}
But then by $f_{2,6}$ we have
\begin{equation}
\beta_2=\alpha_1,\label{eq:beta2=alpha1}
\end{equation}
which contradicts (\ref{eq:gamma4<onehalf}) and
(\ref{eq:beta2=alpha1}). This proves that
\begin{equation}
\gamma_4>\min\{\beta_3,\alpha_2,\gamma_5\}.\label{eq:gamma4nondominant}
\end{equation}
We continue to study the possible dominant values in $f_{2,11}$. There
are two cases to consider.

\begin{itemize}
\item First case : $\gamma_5=\beta_3 <\alpha_2$.\medskip

The possible dominant values are:\smallskip

\begin{eqnarray*}
\ f_{2,6}:\ \ 3\beta_2,\ 2\alpha_1+\gamma_4,\ 2\alpha_1+2\alpha_2\\
\ f_{2,7}:\ 2\beta_2+\beta_3,\ 2\alpha_1+\gamma_5,\
\ \alpha_1+\alpha_2+\gamma_4,\ \alpha_1+3\alpha_2\\
\ f_{2,8}:\ 2\beta_2,\ \beta_2+2\beta_3,\ \alpha_1+\alpha_2+\gamma_5,\
 2\gamma_4,\ 4\alpha_2\\
\ f_{2,9}:\ 3\beta_3,\  \beta_2+\beta_3,\
\ \gamma_5+\gamma_4,\ 2\alpha_2+\gamma_5\\
\ f_{2,10}:\ \beta_2,\ 2\beta_3,\  \gamma_4,\
\ 2\gamma_5\\
\ f_{2,11}:\ \ \beta_3,\ \gamma_5\\
\end{eqnarray*}
\medskip

1) If $\gamma_4\le 2 \gamma_5$ then $\gamma_4<2\alpha_2<\alpha_1$,
\begin{equation}
3\beta_2=2\alpha_1+\gamma_4\label{eq:3beta2=2alpha1gamma4}
\end{equation}
by $f_{2,6}$. Hence
\begin{equation}
\beta_2<\alpha_1,\label{eq:beta2<alpha1}
\end{equation}
so
$$
2\beta_2+\beta_3=\alpha_1+\alpha_2+\gamma_4
$$
by $f_{2,7}$. Thus
\begin{equation}
\beta_3+\alpha_1=\alpha_2+\beta_2.\label{eq:beta3+alpha1=alpha2+beta2}
\end{equation}
By (\ref{eq:beta2<alpha1}) and (\ref{eq:3beta2=2alpha1gamma4}) we have
\begin{equation}
\gamma_4=3\beta_2-2\alpha_1<\beta_2.\label{eq:gamma4<beta2}
\end{equation}
Then $3\beta_3=\gamma_4+\gamma_5$ (that is,
$\gamma_4=2\beta_3$) by $f_{2,9}$ and by $f_{2,8}$ we obtain
$$
\beta_2=\gamma_4,
$$
contradicting (\ref{eq:gamma4<beta2}).
\smallskip

2) Thus $\gamma_4 > 2 \gamma_5$. \medskip

We have $\beta_2> 2\gamma_5$, because otherwise $3\beta_2$ would be the only dominant
value in $f_{2,6}$.  Then the unique dominant value in $f_{2,9}$ is
$3 \beta_3$, a contradiction. This completes the proof in the first case.\medskip

\item Second case: Thus  $\alpha_2\le \gamma_5$ and $\alpha_2\le\beta_3$.\medskip

Then $\gamma_4>\alpha_2$ by (\ref{eq:gamma4nondominant}).\medskip

So the possible dominant values are:\medskip

\begin{eqnarray*}
\ f_{2,6}:\ \ 3\beta_2,\ 2\alpha_1+\gamma_4,\ 2\alpha_1+2\alpha_2\\
\ f_{2,7}:\ 2\beta_2+\beta_3,\ \alpha_1+\alpha_2+\gamma_4,\ \alpha_1+3\alpha_2\\
\ f_{2,8}:\ 2\beta_2,\ \beta_2+2\beta_3,\ 2\gamma_4,\ 4\alpha_2\\
\ f_{2,9}:\ 3\beta_3,\  \beta_2+\beta_3,\ \gamma_5+\gamma_4,\ 2\alpha_2+\gamma_5,\
\gamma_4+\alpha_2,\ 3 \alpha_2\\
\ f_{2,10}:\ \beta_2,\ 2\beta_3,\ \gamma_4,\
2\alpha_2,\ 2\gamma_5\\
\ f_{2,11}:\ \beta_3,\ \alpha_2,\
\gamma_5\\
\end{eqnarray*}
1) If
\begin{equation}
\gamma_4<2\alpha_2\label{eq:gamma4<alpha2bis}
\end{equation}
then
\begin{equation}
\gamma_4=\beta_2\label{eq:gamma4=beta2bis}
\end{equation}
by $f_{2,10}$. From $f_{2,6}$ we obtain the equality (\ref{eq:3beta2=2alpha1gamma4}), which implies
\begin{equation}
\beta_2=\alpha_1,
\end{equation}
which contradicts (\ref{eq:gamma4<alpha2bis}) and (\ref{eq:gamma4=beta2bis}).\smallskip

2) Suppose $\gamma_4>2\alpha_2$. By looking at the dominant terms of $f_{2,6}$ and $f_{2,7}$ we obtain again the equality (\ref{eq:beta3+alpha1=alpha2+beta2}). If $\beta_2\le2\alpha_2$ then by $f_{2,6}$ we would have $\alpha_2\ge\frac{1}{2}\alpha_1$, which is false. Hence $\beta_2>2\alpha_2$. Then by $f_{2,6}$ we have $\beta_2<\alpha_1$ and now (\ref{eq:beta3+alpha1=alpha2+beta2}) implies $\beta_3>\alpha_2$. Then the only possible dominant value in $f_{2,8}$ is $4\alpha_2$, a contradiction.\smallskip

3) So $\gamma_4=2\alpha_2$.\medskip

Then
\begin{equation}
\beta_2>2\alpha_2\label{eq:beta2>2alpha2}
\end{equation}
(if not  $3\beta_2$ would be the only dominant value in $f_{2,6}$). Using $f_{2,6}$ and $f_{2,7}$ we
see that
\begin{equation}
2\beta_2+\beta_3\ge\frac{4}{3}(\alpha_1+\alpha_2)+\alpha_2=
\alpha_1+2\alpha_2+\frac{1}{3}(\alpha_1+\alpha_2)>\alpha_1+3\alpha_2.\label{eq:2beta2+beta3vsalpha1+3alpha2}
\end{equation}
Let us do another trichotomy:

$-$ A) Suppose $\alpha_2=\gamma_5 <\beta_3$
\medskip

The possible dominant values are:
\begin{eqnarray*}
\ f_{2,6}:\ 3\beta_2 \ge 2\alpha_1+\gamma_4=2\alpha_1+2\alpha_2\\
\ f_{2,7}: \alpha_1+\alpha_2+\gamma_4= \alpha_1+3\alpha_2\\
\ f_{2,8}:\ 2\gamma_4=\ 4\alpha_2\\
\ f_{2,9}:\ \gamma_5+\gamma_4=\ 2\alpha_2+\gamma_5=
\gamma_4+\alpha_2=\ 3 \alpha_2\\
\ f_{2,10}:\  \gamma_4=
2\alpha_2=2\gamma_5\\
\ f_{2,11}:\ \alpha_2=
\gamma_5\\
\end{eqnarray*}
If $3\beta_2> 2\alpha_1+2\alpha_2$ then the first equations of any
wedge with $b_4\not=0$ are:
\begin{eqnarray*}
\  &\bar f_{2,6}=2i a_1^2( c_4-i a_2^2)=0\\
\ &\bar f_{2,7}=4i a_1 a_2( c_4-i a_2^2)=0\\
&\bar f_{2,8}=( c_4+i a_2^2)( c_4-i a_2^2)=0\\
&\bar f_{2,9}=\ ( c_4+i a_2^2)( c_5-2i a_2 a_3)+( c_5+2i a_2 a_3)( c_4-i
a_2^2)=0\\
&\bar f_{2,10}=( c_4+i a_2^2)( c_6-i a_3^2)+( c_5+2i a_2 a_3)( c_5-2i a_2
a_3)+( c_6+i a_3^2)( c_4-i a_2^2)=0\\
&\bar f_{2,11}=( c_5+2i a_2 a_3)( c_6-i a_3^2)+( c_6+i a_3^2)( c_5-2i a_2
a_3)=0.\\
& f_{6,12}=a_3^4+b_4^3+c_6^2=0.
\end{eqnarray*}

Thus $c_4-i a_2^2=0$ and the last four equations become:

 \begin{eqnarray*}
&\bar f_{2,9}=\ ( c_4+i a_2^2)( c_5-2i a_2 a_3)=0\\
&\bar f_{2,10}=( c_4+i a_2^2)( c_6-i a_3^2)+( c_5+2i a_2 a_3)( c_5-2i a_2
a_3))=0\\
&\bar f_{2,11}=( c_5+2i a_2 a_3)( c_6-i a_3^2)+( c_6+i a_3^2)( c_5-2i a_2
a_3)=0.\\
& f_{6,12}=a_3^4+b_4^3+c_6^2=0
\end{eqnarray*}
Since $b_4\ne0$, we have $c_6-i a_3^2\ne0$ and $c_6-i a_3^2\ne0$. As
well, $c_5\ne0$, $a_2\ne0$, $c_4\ne0$ which is incompatible with the
above equations.\smallskip

Thus $3\beta_2= 2\alpha_1+2\alpha_2$. The first equations of the wedge are: 
 \begin{eqnarray*}
\  &\bar f_{2,6}= b_2^3+2i a_1^2( c_4-i a_2^2)=0\\
\ &\bar f_{2,7}=4i a_1 a_2( c_4-i a_2^2)=0\\
\end{eqnarray*}
Thus as $a_1\ne0$ and $a_2\ne0$, we have $ c_4-i a_2^2=0$ and then
$b_2=0$ (not allowed by definition).\\
\smallskip

Thus the case A) is impossible and $\alpha_2=\beta_3$.
\medskip

$-$  B) $\alpha_2=\beta_3<\gamma_5$.
\medskip

Using (\ref{eq:beta2>2alpha2}) and
(\ref{eq:2beta2+beta3vsalpha1+3alpha2}), we see that the possible
dominant values are:
\begin{eqnarray*}
\ f_{2,6}:\ 3\beta_2 \ge 2\alpha_1+\gamma_4=2\alpha_1+2\alpha_2\\
\ f_{2,7}: \alpha_1+\alpha_2+\gamma_4= \alpha_1+3\alpha_2\\
\ f_{2,8}:\ 2\gamma_4=\ 4\alpha_2\\
\ f_{2,9}:\ 3\beta_3=
\gamma_4+\alpha_2=\ 3 \alpha_2\\
\ f_{2,10}:\  \gamma_4=
2\alpha_2=2\beta_3\\
\ f_{2,11}:\ \alpha_2=\beta_3\\
\end{eqnarray*}
Suppose that $3\beta_2 = 2\alpha_1+\gamma_4=2\alpha_1+2\alpha_2$, then
the first equations of the wedge are:

 \begin{eqnarray}
\  &\bar f_{2,6}= b_2^3+2i a_1^2( c_4-i a_2^2)=0\label{eq:barf26}\\
\ &\bar f_{2,7}=4i a_1 a_2( c_4-i a_2^2)=0\label{eq:barf27}
\end{eqnarray}

Since $a_1\ne0$ and $a_2\ne0$, we obtain $c_4-i a_2^2=b_2=0$, which
gives the desired contradiction.\smallskip

Therefore $3\beta_2 >2\alpha_1+\gamma_4=2\alpha_1+2\alpha_2$. The
first equations are:

\begin{eqnarray*}
\  &\bar f_{2,6}= 2i a_1^2( c_4-i a_2^2)=0\\
\ &\bar f_{2,7}=4i a_1 a_2( c_4-i a_2^2)=0\\
&\bar f_{2,8}=( c_4+i a_2^2)( c_4-i a_2^2)=0\\
&\bar f_{2,9}=  b_3^3+ 4a_2^3a_3=0\\
&\bar f_{2,10}=( c_4+i a_2^2)( c_6-i a_3^2)+( c_4-i a_2^2)( c_6+i
a_3^2)+3b_3^2b_4=0\\
&\bar f_{2,11}=4a_2a_3^3+3b_3b_4^2=0.\\
& f_{6,12}=a_3^4+b_4^3+c_6^2=0
\end{eqnarray*}
First we see that $c_4-i a_2^2=0$ and the equations become:
\begin{eqnarray*}
& c_4-i a_2^2=0\\
&\bar f_{2,9}=  b_3^3+ 4a_2^3a_3=0\\
&\bar f_{2,10}=2i a_2^2( c_6-i a_3^2)+3b_3^2b_4=0\\
&\bar f_{2,11}=4a_2a_3^3+3b_3b_4^2=0.\\
& f_{6,12}=a_3^4+b_4^3+c_6^2=0
\end{eqnarray*}

By XMaple (one can also do it by hand, as in $f_{2,11}$ the equation is
linear in $b_3$ and $a_2$), these equations imply that $b_3=0$, which is not
allowed by definition of $b_3$. Thus case B) is also
impossible and the only remaining case to consider is
\medskip

$-$ C)
$$
\alpha_2=\beta_3=\gamma_5<\frac{1}{2}\alpha_1,\quad \frac23\alpha_1\le\beta_2\ \text{and}
\ \gamma_4=2\alpha_2.
$$
Using (\ref{eq:beta2>2alpha2}) and (\ref{eq:2beta2+beta3vsalpha1+3alpha2}), we see that the possible dominant values are:
\begin{eqnarray*}
\ f_{2,6}:\ 3\beta_2 \ge 2\alpha_1+\gamma_4=2\alpha_1+2\alpha_2\\
\ f_{2,7}: \alpha_1+\alpha_2+\gamma_4= \alpha_1+3\alpha_2\\
\ f_{2,8}:\ 2\gamma_4=\ 4\alpha_2\\
\ f_{2,9}:\ 3\beta_3=
\gamma_4+\alpha_2=\ 3 \alpha_2=\gamma_4+\gamma_5=2\alpha_2+\gamma_5\\
\ f_{2,10}:\  \gamma_4=
2\alpha_2=2\beta_3=2\gamma_5\\
\ f_{2,11}:\ \alpha_2=
\beta_3=\gamma_5\\
\end{eqnarray*}

If $3\beta_2=2\alpha_1+2\alpha_2$ then the first two equations of the wedge are (\ref{eq:barf26}) and (\ref{eq:barf27}). We obtain the same contradiction as before: since $a_1\ne0$ and $a_2\ne0$, we have $ c_4-i a_2^2=0$ and then
$b_2=0$ (not allowed by definition).\medskip

Finally, it remains to solve the case when
$$
3\beta_2 > 2\alpha_1+2\alpha_2.
$$
The equations of the wedge are:

 \begin{eqnarray*}
\  &\bar f_{2,6}= 2i a_1^2( c_4-i a_2^2)=0\\
\ &\bar f_{2,7}=4i a_1 a_2( c_4-i a_2^2)=0\\
&\bar f_{2,8}=( c_4+i a_2^2)( c_4-i a_2^2)=0\\
&\bar f_{2,9}=\ ( c_4+i a_2^2)( c_5-2i a_2 a_3)+b_3^3=0\\
&\bar f_{2,10}=( c_4+i a_2^2)( c_6-i a_3^2)+( c_5+2i a_2 a_3)( c_5-2i a_2
a_3)+3b_3^2b_4=0\\
&\bar f_{2,11}=( c_5+2i a_2 a_3)( c_6-i a_3^2)+( c_6+i a_3^2)( c_5-2i a_2
a_3)+3b_3b_4^2=0.\\
& f_{6,12}=a_3^4+b_4^3+c_6^2=0
\end{eqnarray*}

As by definition $a_1 \not = 0$, we have $c_4-i a_2^2=0$ and the last four equations become:
 \begin{eqnarray*}
&\bar f_{2,9}=\ 2i a_2^2( c_5-2i a_2 a_3)+b_3^3=0\\
&\bar f_{2,10}=2i a_2^2( c_6-i a_3^2)+( c_5+2i a_2 a_3)( c_5-2i a_2
a_3)+3b_3^2b_4=0\\
&\bar f_{2,11}=( c_5+2i a_2 a_3)( c_6-i a_3^2)+( c_6+i a_3^2)( c_5-2i a_2
a_3)+3b_3b_4^2=0.\\
& f_{6,12}=a_3^4+b_4^3+c_6^2=0
\end{eqnarray*}

By XMaple, one obtains $b_3=0$, which is not allowed by definition of $b_3$. So in this last case one cannot construct the wedge either. This completes the proof of the last non-inclusion. $\Box$

\end{itemize}

\end{itemize}

\end{enumerate}

\end{document}